\documentclass{article}
\usepackage{lipsum}
\usepackage[utf8]{inputenc}
\usepackage{fullpage}
\usepackage{amsmath}
\usepackage{amssymb}
\usepackage{amsthm}
\usepackage{amscd}
\usepackage{amsfonts}
\usepackage{graphicx}
\usepackage{fancyhdr}
\usepackage{authblk}
\usepackage{mathrsfs}
\usepackage{comment}
\title{Prasad's Conjecture about dualizing involutions}
\author{Prashant Arote and  Manish Mishra}
\date{}

\AtEndDocument{\bigskip{\footnotesize%
  \textsc{Department  of Mathematics, Indian Institute of Science Education and Research Pune,}\par
  \textsc{Dr. Homi Bhabha Road, Pune 411008, India} \par  
  P.Arote: \texttt{prashantarote123@gmail.com} \par
  M. Mishra: \texttt{manish@iiserpune.ac.in} 
}}

\theoremstyle{plain} \numberwithin{equation}{section}
\newtheorem{theorem}{Theorem}[section]
\newtheorem{corollary}[theorem]{Corollary}
\newtheorem{lemma}[theorem]{Lemma}
\newtheorem{proposition}[theorem]{Proposition}
\theoremstyle{definition}
\newtheorem{definition}[theorem]{Definition}

\newtheorem{remark}[theorem]{Remark}

\newtheorem{example}[theorem]{Example}
\newtheorem{Fact}[theorem]{Fact}
\newtheorem{notation}[theorem]{Notation}

\usepackage[colorlinks=true]{hyperref}

\makeatletter
\@ifpackageloaded{hyperref}%
  {\newcommand{\mylabel}[2]
    {\protected@write\@auxout{}{\string\newlabel{#1}{{#2}{\thepage}%
      {\@currentlabelname}{\@currentHref}{}}}}}%
  {\newcommand{\mylabel}[2]
    {\protected@write\@auxout{}{\string\newlabel{#1}{{#2}{\thepage}}}}}
\makeatother

\usepackage{amsmath}
\usepackage{tikz-cd}
\usepackage[all,cmtip]{xy}
\usepackage{graphicx}
\usepackage{amsmath}
\usepackage{wasysym}
\usepackage{amssymb}
\usepackage{amsthm}
\usepackage{commath}
\usepackage{mathtools}
\usepackage{regexpatch,fancyvrb,xparse}
\usepackage{tikz-cd}
\usepackage[utf8]{inputenc}
\usepackage[nottoc]{tocbibind}
\usepackage{rotating}
\usepackage[margin=1.0in]{geometry}
\usepackage{csquotes}

\usepackage{lipsum}

\newcommand\blfootnote[1]{%
  \begingroup
  \renewcommand\thefootnote{}\footnote{#1}%
  \addtocounter{footnote}{-1}%
  \endgroup
}

\newcommand{\Ind}{\operatorname{Ind}}

\newcommand{\Ql}{\mathbb{Q}_{\ell}}

\newcommand{\Qlcl}{\overline{\mathbb{Q}}_{\ell}}

\newcommand{\E}{\mathcal{E}}
\newcommand{\Hom}{\operatorname{Hom}}

\newcommand{\Irr}{\operatorname{Irr}}

\newcommand{\Aut}{\operatorname{Aut}}

\newcommand{\fq}{\mathbb{F}_{q}}
\newcommand{\Id}{\operatorname{Id}}
\newcommand{\Uch}{\operatorname{Uch}}
\newcommand{\J}{\mathbb{J}}

\newcommand{\p}{\mathcal{P}}
\newcommand{\Stab}{\operatorname{Stab}}
\newcommand{\ad}{\operatorname{ad}}
\newcommand{\Iso}{\operatorname{Iso}}
\newcommand{\GL}{\operatorname{GL}}
\newcommand{\Sp}{\operatorname{Sp}}
\newcommand{\U}{\operatorname{U}}
\newcommand{\SO}{\operatorname{SO}}
\newcommand{\SL}{\operatorname{SL}}

\newcommand{\der}{\operatorname{der}}

\newcommand{\AS}{\widehat{A(s)^{F^*}}}
\newcommand{\Gad}{G_{\ad}^F/G^F}

\begin{document}


\maketitle

\begin{abstract} 
Let $G$ be a connected reductive group defined over a finite field $\fq$ with corresponding Frobenius $F$.
Let $\iota_G$ denote the duality involution defined by D. Prasad under the hypothesis $2\mathrm{H}^1(F,Z(G))=0$, where $Z(G)$ denotes the center of $G$.  
We show that for each irreducible character $\rho$ of $G^F$, the involution $\iota_G$ takes $\rho$ to its dual $\rho^{\vee}$ if and only if for a suitable Jordan decomposition of characters, an associated unipotent character $u_\rho$ has Frobenius eigenvalues $\pm$ 1. As a corollary, we obtain that if $G$ has no exceptional factors and satisfies $2\mathrm{H}^1(F,Z(G))=0$, then the duality involution $\iota_G$ takes $\rho$ to its dual $\rho^{\vee}$ for each irreducible character $\rho$ of $G^F$. Our results resolve a finite group counterpart of a conjecture of D.~Prasad.
\end{abstract}

\section{Introduction}
\blfootnote{The first named author was suported by NBHM. The Second named author was supported by SERB Core Research Grant CGR/2022/000415}
A well known result due to Gelfand and Kazhdan for general linear groups states that the involution - transpose inverse - takes every irreducible representation to its dual representation \cite{GelKazh72}. The purpose of this paper is to produce a generalization of this result for finite reductive groups. 

Let $G$ be a connected reductive group defined over a field $K$.
If $G$ is quasi-split, then given a $K$-pinning $\mathcal{P}$ of $G$, D. Prasad constructed an involutive
automorphism $\iota_{G,\mathcal{P}}$ of $G(K)$ generalizing an involution
used by Moeglin, Vign\'eras and Waldspurger for classical groups
over local fields \cite[Chap. IV, \S II]{MVW87}. The latter has the property that it takes any irreducible
admissible representation $\pi$ of $G(K)$, when $G$ is classical and $K$ is local, to its contragradient
$\pi^{\vee}$. The result in loc.~cit. was reproved and slightly extended in \cite{RoVin18} using different methods. 

The involution  $\iota_{G,\p}$ defined by Prasad is independent
of the pinning $\mathcal{P}$ if $2\mathrm{H}^{1}(K,Z(G))=0$, where
$Z(G)$ denotes the center of $G$. When $K$ is a finite field, this
is also a necessary condition. From now on, $K$ is a finite field $\fq$ of order $q$ and $F$ is the corresponding Frobenius on $G$. Recall here that finite reductive groups are necessarily quasi-split. In this paper, we show that under the hypothesis $2\mathrm{H}^{1}(F,Z(G))=0$, the duality involution takes each irreducible representation of $G^F$ to its dual representation if $G$ contains no factors of exceptional type. More generally, under the above hypothesis on $Z(G)$, we translate the condition $\rho\circ\iota_G=\rho^{\vee}$ for each irreducible character $\rho$ of $G^F$ to an equivalent condition on the unipotent counterpart of $\rho$ for a suitable Jordan decomposition of characters. The condition on unipotent counterpart (namely the $\emph{Frobenius eigenvalue}$ being $\pm$ 1) is easily computable and can be looked up from lists given in the literature. Our results resolve a finite group counterpart of a conjecture of D. Prasad \cite[Conjecture 1]{dualityinvolution} (see Remark \ref{rem:prasad} for comparison with the $p$-adic case).

To discuss our results further, let us first recall \cite{Jordandecopo} the notion of Lusztig's Jordan decomposition of characters of $G^F$. To each conjugacy class of a semisimple element $s$ in the dual
group ${G^*}^{F^*}$ of $G$, Lusztig associated a subset $\mathcal{E}(G^{F},s)$
of the set $\mathrm{Irr}(G^{F})$ of irreducible complex characters
of $G^{F}$. The subset $\mathcal{E}(G^{F},s)$ is called the \textit{Lusztig
series} associated to $s$. The collection of all Lusztig series of
$G$ partition the set $\mathrm{Irr}(G^{F})$. The elements of $\mathcal{E}(G^{F},1)$
are called the \textit{unipotent characters} of $G$. We will denote
them by $\mathrm{Uch}(G^{F})$. When $G$ has connected centre, Lusztig
established that there exists a bijection $J_{s}^{G}$ between the
sets $\mathcal{E}(G^{F},s)$ and $\mathrm{Uch}(C_{G^{*}}(s)^{F^{*}})$
satisfying a certain natural property. Here $C_{G^{*}}(s)$ denotes
the centralizer of $s$ in $G^{*}$. 
The collection of bijections
$J_{s}^{G}$ thus defines the notion of \textit{Jordan decomposition}
of characters. 
This notion can be extended to include groups with
disconnected centre. 
In the case of  groups with disconnected centre, $J_{s}^{G}$
is a surjection from $\mathcal{E}(G^F,s)$ onto equivalence classes
of $\mathrm{Uch}(C_{G^{*}}(s)_{\circ}^{F^*})$, where equivalence is given
by the action of the group $C_{G^{*}}(s)^{F^{*}}/C_{G^{*}}(s)_{\circ}^{F^{*}}$.
Here $(-)_{\circ}$ denotes the identity component of $(-)$. 

Now fix a pinning $\p$ of $G$ and let $G^*$ be the dual group determined by this pinning with associated Frobenius $F^*$. Let $\sigma$ be an $F$-automorphism of $G$ which respects the pinning $\p$ and let $\sigma^*$ be the associated $F^*$-automorphism of $G^*$. Our first result (see Theorem \ref{thm:connectedcentre}) proves the compatibility of Jordan decomposition with the corresponding actions of $\sigma$ and $\sigma^*$ when $G$ has connected center. More precisely: 
\begin{theorem}
A Jordan decomposition can be chosen such that for each $\rho\in \E(G^F,s)$, $$J_{s^{-1}}(\rho^{\vee})=J_s(\rho)^{\vee}.$$ Moreover, if $Z(G)$ is connected, then,
$$J_{{\sigma^*}^{-1}(s)}(\rho\circ{\sigma^{-1}})=J_s(\rho)\circ{\sigma^*}^{-1}.$$
\end{theorem}

For a character $\rho$, if $\rho\circ\iota_{G,\p}$ is independent of $\p$, we will drop the subscript $\p$ for simplicity of notation even if $\iota_{G,\p}$ itself is dependent on $\p$. Using the above Theorem, we show:

\begin{theorem}  Let $G$ be a connected reductive group defined over $\fq$ satisfying
 $2\mathrm{H}^1(F,Z(G))=0$. Then a Jordan decomposition can be chosen such that for each $\rho\in \E(G^F,s)$,
$$J_{s^{-1}}(\rho\circ{\iota_G})=J_s(\rho)\circ{\iota_{G^*}}.$$
\end{theorem}
We now state our main result (see Theorems \ref{thm:iotadual} and \ref{thm:main}). 
\begin{theorem}\label{thm:intro_main}  Let $G$ be a connected reductive group defined over $\fq$ satisfying
     $2\mathrm{H}^1(F,Z(G))=0$. Then a Jordan decomposition can be chosen such that for each irreducible character $\rho$ of $G^F$,
    $$\rho\circ\iota_{G}=\rho^{\vee}\ \ \mbox{if and only if}\ \ u_\rho\circ\iota_{G^*}=u_\rho^{\vee},$$ for any, and therefore all $u_\rho\in J_s(\rho),$ $$\mbox{~if  and only if the Frobenius eigenvalues of~} u_\rho \mbox{~are~} \pm 1,$$ for any, and therefore all $u_\rho\in J_s(\rho)$.
\end{theorem}
As an immediate corollary, we obtain (see Corollary \ref{thm:classical}),
\begin{theorem}
    Let $G$ be a connected reductive group defined over $\fq$ satisfying $2\mathrm{H}^1(F,Z(G))=0$. Assume further that $G$ contains no factors of exceptional type. Then for all irreducible representation of $G^F$, we have, $$\rho\circ\iota_G=\rho^{\vee}.$$
\end{theorem}

The main ingredient of our proofs is the Theorem of Digne and Michael about uniqueness of Jordan decomposition (see Theorem \ref{Jordandecompositionconnectedcentre}). 

For irreducible generic representations (i.e., those representations which appear in some Gelfand-Graev representation), the hypothesis on $Z(G)$ can be dropped. The following result (see Theorem \ref{thm:generic}) follows quite easily and independently without requiring any of the previously stated theorems.

\begin{theorem}\label{thm:intro_generic}
    Let $G$ be any connected reductive group over $\fq$. Let $\p$ be the pinning of $G$ corresponding to a given algebraic Whittaker datum $(B,\psi)$, where $\psi$ is a non-degenerate character of the unipotent radical of the given $F$-stable Borel $B$. Then for a $\psi$-generic representation $\rho$ of $G^F$, $$\rho\circ\iota_{G,\p}=\rho^{\vee}.$$
\end{theorem}

Theorem \ref{thm:intro_generic} was proved by Chang Yang \cite[Theorem 1.3]{Changyang} using the theory of Bessel functions. The result in loc.~cit. covers groups over non-archimedean local fields as well. 

Theorem \ref{thm:intro_main} allows us to associate a sign in a natural way to most irreducible character of $G^F$ under the hypothesis of the same Theorem. This is observed in Section \ref{sec:sign}.
\subsection*{Notation}\label{sec:notation}

Given a homomorphism of finite groups $f:G\rightarrow G'$, if $f(G)$ is normal in $G'$, we will often abuse notation and write $G'/G$ for $G'/f(G)$ whenever there is no possibility of confusion. If $\rho$ is a representation of $G'$, we will abuse notation to write $\rho|_{G}$ to denote $\rho|_{f(G)}$.

\section{Review}
In this section, we recall briefly the notions of the \emph{Duality involution} and the \emph{Jordan decomposition of characters} of finite reductive groups. 
We also recall the notion of \textit{dual morphism}.
For more details we refer to \cite{dualityinvolution}, 
 \cite{Jordandecopo}, \cite{Jordandecodisconnectedcase}, 
 \cite{Book:DigneandMichel}, \cite{Book:GeckandMalle}, \cite{actionofauto}.

\subsection{Duality of automorphism}
Let $G$ be a connected reductive group defined over a finite field $\fq$ and let $F:G\rightarrow G$ be the corresponding Frobenius morphism.
Fix an $F$-stable Borel subgroup $B_{0}$ of $G$ and a maximal $F$-stable torus $T_{0}\leq B_{0}$.
With respect to $T_0$, we have the root datum of $G$,
$(X(T_0), \Phi, X^{\vee}(T_0), \Phi^{\vee})$ of $G$.
Here $X(T_0)$, $X^{\vee}(T_{0})$ denotes the character group and the cocharacter group of $T_0$.
If $\Phi\subset X(T_0)$ are the roots of $G$ with
respect to $T_0$, then associated to each $\alpha\in\Phi$, we have a corresponding $1$-dimensional root subgroup $U_{\alpha}\leq G$.
The set of all roots $\alpha\in\Phi$ such that $U_{\alpha}\subset B_{0}$ is a positive system of roots which contains a unique simple system of roots, say $\Delta$.
Then $\mathcal{R}(G)=(X(T_0), \Phi, \Delta, X^{\vee}(T_0), \Phi^{\vee}, \Delta^{\vee})$ a based root datum determined by the triple $(G, B_0, T_0)$.
By the existence theorem for connected reductive groups, there exists a triple $(G^*,B_0^*,T_0^*)$, with a based root datum $\mathcal{R}(G^*)$, such that we
have an isomorphism $\mathcal{R}(G^*)\rightarrow\mathcal{R}(G)^{\vee}$ of based root data, where $\mathcal{R}(G)^{\vee}$ denotes the dual root datum of $\mathcal{R}(G)$. 
We say that the pair $(G, G^*)$ is a dual pair
or that $G^*$ is dual to $G$.
\begin{notation}
Let $G$ be a connected reductive group defined over $\fq$ and let $(G,B_0,T_0)$ be a triple as above.
Then $\Iso(G, B_0, T_0)$ denotes the set of isogenies  $\phi:G\rightarrow G$  such that $\phi(B_0)=B_0$ and $\phi(T_0)=T_0$ and $\Iso((G, B_0, T_0),F)$ denotes the subset of $\Iso(G, B_0, T_0)$ of those automorphisms which commute with $F$.
We will denote by $\Aut(G, B_0, T_0)$, the set of all automorphisms $\phi:G\rightarrow G$ such that $\phi(B_0)=B_0$ and $\phi(T_0)=T_0$.
In addition, we will denote by $\Aut((G, B_0, T_0),F)\subseteq\Aut(G, B_0, T_0)$, those automorphisms that commute with $F$. 
\end{notation}
\begin{theorem}[Chevalley]\label{thm chevally}
Let $G$ be a connected reductive group defined over $\fq$ and let $(G,B_0,T_0)$ be a triple as above.
Let $\mathcal{R}(G)$ be a based root datum of $G$ determined by the triple $(G,B_0,T_0)$.
Then there is a bijection between $\Iso(G,B_0,T_0)/T_0$ and $\Iso_p(\mathcal{R}(G))$ and there is a surjection $\Iso((G,B_0,T_0),F)\rightarrow\Iso_p(\mathcal{R}(G),F)$, where $\Iso_p(\mathcal{R}(G))$ (resp. $\Iso_p(\mathcal{R}(G),F)$ denotes the set of $p$-isogenies (resp. $p$-isogenies which commute with $F$) of the root datum $\mathcal{R}(G)$.
\end{theorem}
\begin{Fact}[{\cite[\textsection 5]{actionofauto}}]
Let $(G,G^*)$ be a dual pair.
Then there is a natural bijection,
$$\Iso(G,B_0,T_0)/T_0\rightarrow\Iso(G^*,B_0^*,T_0^*)/T_0^*,$$
between the orbits.
An isogeny $\sigma\in\Iso(G,B_0,T_0)$ is said to be \textit{dual} to $\sigma^*\in\Iso(G^*,B_0^*,T_0^*)$ if their orbits corresponds under the above bijection. 
\end{Fact}
\begin{lemma}[{\cite[Lemma 5.5]{actionofauto}}]\label{lemmaexistanceofdualauto}
    If $\sigma^*\in\Iso(G^*,B_0^*,T_0^*)$ is an isogeny dual to $\sigma\in\Iso(G,B_0,T_0)$, then $\sigma$ is injective if
and only if $\sigma^*$ is injective.
Furthermore, if $\sigma$ is injective, and hence bijective, then $\sigma^{*^{-1}}$ is dual to $\sigma^{-1}$.
\end{lemma}
\begin{definition}
Let $\sigma\in\Aut((G,B_0,T_0),F)$.
Then by Theorem \ref{thm chevally}, there exists an automorphism $\sigma^{*}\in\Aut((G^*,B_0^*,T_0^*),F^*)$ whose orbit is in bijection with the orbit of $\sigma$.
We say that $\sigma^*$ is \textit{dual} to $\sigma$.
\end{definition}

\subsection{Duality Involution $\iota_{G,\p}$}
In this Section, we  recall  Prasad's construction of the  duality involution $\iota_{G,\p}$, for any connected reductive group $G$ defined over $\fq$ with a fixed pinning $\p$.
For more details, we refer to \cite{dualityinvolution}.

\begin{definition}[Chevalley involution]
Fix a pinning $\mathcal{P}=(G,B,T,\{X_{\alpha}\})$ of $G$.
Then there exists an involution $c_{G,\p}$ of $G$ defined on $\fq$ that fixes a pinning $\p$ and acts as $t\mapsto w_G(t)^{-1}$ on $T$, where $w_G$ is the longest element in the Weyl group of $T$ taking $B$ to $\Bar{B}$, the opposite Borel subgroup.
This involution is called the Chevalley involution associated to a pinning $\p$ on $G$.
\end{definition}
\begin{remark}[{\cite[\textsection 3]{dualityinvolution}}]
This involution induces an isomorphism $-w_{G}$ on the based root datum 
 of $G$ associated a pinning $\p$.
In addition, one can observe that $c_{G,\p}=1$ if and only if $-1\in W_{G}(T)$ . 
Here $-1$ denotes the negative identity map on the root datum, $\emph{i.e.},\ \alpha\mapsto-\alpha$.
Therefore, for all simple groups except those of type $A_n\ (n\geq 2)$, $D_{2n+1}$, $E_6$, $c_{G,\p} = 1$.
\end{remark}
\begin{Fact}\label{rmk:dual of duality}
Let $(G,{G}^*)$ be a dual pair.
Fix a pinning $\mathcal{P}=(G,B,T,\{X_{\alpha}\})$ on $G$. This determines a pinning $\p^*$ on $G^*$.
Then
$c_{G,\p}\in\Aut((G,B,T),F)$ and $c_{G,\p}$ is dual to $c_{G^*,\p^*}\in\Aut((G^*,B^*,T^*),F^*)$.
\end{Fact}
\begin{Fact}[{\cite[\textsection 3]{dualityinvolution}}]
Given a pinning $\mathcal{P}=(G,B,T,\{X_{\alpha}\})$ of $G$, there exists a unique element $t\in T_{\ad}^F$ such that for any 
simple root $\alpha$ determined by $\p$, $\alpha(t)=-1$.
Let $\iota_{-}\coloneqq \ad(t)$, the inner automorphism of $G$ induced by $t$.
Then this automorphism $\iota_{-}$ commutes with a Chevalley involution $c_{G,\p}$.
\end{Fact}
\begin{definition}[Duality Involution {\cite[Def.1]{dualityinvolution}}]
Let $G$ be a connected reductive group defined over $\fq$ and let $\mathcal{P}=(G,B,T,\{X_{\alpha}\})$ be a pinning of $G$. 
Define $\iota_{G,\p}\coloneqq \iota\circ c_{G,\p}$, where $c_{G,\p}$ and $\iota_{-}$ are constructed as above.
The involution $\iota_{G,\p}$ is called the \textit{duality involution}.
\end{definition}
\begin{example}
Let $G=\GL_n$ then one can check that duality involution is $g\mapsto (g^{t})^{-1}$ transpose-inverse morphism.
\end{example}
Let $\p,\p'$ be two pinnings of $G$.
Since $G_{\ad}^F$ acts transitively on the set of pinnings of $G$, there exists $g\in G_{\ad}^F$ such that $\p'=\ad(g)(\p)$.
The corresponding duality involutions are related by a following relation:
$$\iota_{G,\p'}=\ad(g)\circ\iota_{G,\p}\circ \ad(g)^{-1}.$$
\begin{proposition}[{\cite[Prop.1]{dualityinvolution}}] \label{prop pinnindependence}
Let $G$ be a connected reductive group defined over $\fq$.
The automorphism $\iota_{G,\p}$ is independent of the choice of $\p$ (as an element of $\Aut(G,F)/G^F$) if and only if  $2\mathrm{H}^1(F,Z(G))=0$.
\end{proposition}

\begin{remark}
    If $Z(G)$ satisfies any of the conditions of \cite[Corollary 1]{dualityinvolution}, then duality involution is independent of the choice of the pinning. Equivalently, $2H^1(F,Z(G))=0$. Thus for all classical groups $G=\GL_n, \U_n, \SO_n, \Sp_{2n}$, the duality involution is independent of the choice of a pinning. On the other hand, for $\SL_n$, it does depend on the choice of a pinning.
\end{remark}
\begin{notation}
If the condition in Prop. \ref{prop pinnindependence} is satisfied, then we will denote the duality involution by $\iota_{G}$.
\end{notation}

\begin{definition}[Algebraic Whittaker datum]
An algebraic Whittaker datum is a pair $(B,\psi)$, where $B=TU$ is an $F$-stabl Borel subgroup of $G$ with unipotent radical $U$ and $\psi:U^F\rightarrow \fq$ is a non-degenerate character of $U^F$, i.e., the restriction of $\psi$ to any
simple root subgroup  $U_{\alpha}^F\subset U$ is nontrivial.
\end{definition}
\begin{remark}[{\cite[Remark 1]{dualityinvolution}}]
There is a bijection between the set of algebraic Whittaker datum on $G$ and $U^F$-conjugacy classes of pinnings $(G,B,T,\{X_{\alpha}\})$. This bijection is obtained by sending $(B,\psi)$ to the collection of elements $\{X_{\alpha}\}$, where $\psi(X_{\alpha})\neq1$ for all simple roots $\alpha$.
\end{remark}
In  view of the above remark, there is an another construction of the involution $\iota_{G,\p}$ using the algebraic Whittaker datum corresponding to the pinning $\p$.
Suppose, $\psi:U^F\rightarrow \fq$ is an algebraic Whittaker datum.
Then the group $G$ has an automorphism $\iota_{B,T,U,\psi}$ which commutes with $F$ and it has the following properties:
\begin{enumerate}
    \item $\iota_{B,T,U,\psi}$ maps the pair $(T,B)$ to itself.
    \item $\psi\circ\iota_{B,T,U,\psi}=\psi^{-1}$.
    \item If $\alpha$ is a simple root of $T$ on $U$ then $\iota_{B,T,U,\psi}(\alpha)=-w_G(\alpha)$, where $w_G$ is the longest element in the Weyl group of $G$.
\end{enumerate}
If $(B,\psi)$ corresponds to pinning $\p$ then one can take $\iota_{G,\p}=\iota_{B,T,U,\psi}$.
\begin{definition}[Gelfand-Graev representation]
Let $(B,\psi)$ be an algebraic Whittaker datum. 
The \emph{Gelfand-Graev representation} corresponding to $\psi$ is denoted by $\Gamma_{\psi}$ and defined as $\Gamma_{\psi}=\Ind_{U^F}^{G^F}(\psi)$.
 The irreducible constituents of $\Gamma_{\psi}$ are called  $\psi$-\textit{generic} representations.
\end{definition}
\begin{definition}[semisimple character] 
 Let $\rho$ be an irreducible representation of $G^F$. 
 Then $\rho$ is said to be  a semisimple character of $G^F$ if $\pm D_G(\rho)$ is a $\psi$-generic  representation for some $\psi$, where $D_G$ is the Alvis–Curtis–Kawanaka–Lusztig duality operator on the space of class functions on $G^F$.
\end{definition}
\begin{lemma}\label{lemma generic}
Let $(B,\psi)$ be an algebraic Whittaker datum and $\p$ be the corresponding pinning.
Then $\Gamma_{\psi}\circ\iota_{G,\p}=\Gamma_{\psi^{-1}}$.
As a consequence, if $\rho$ is a $\psi$-generic representation, then $\rho\circ\iota_{G,\p}$ is a $\psi^{-1}$-generic representation.
\end{lemma}
\begin{proof}
    Consider, $\Gamma_{\psi}\circ\iota_{G,\p}=\Gamma_{\psi}\circ\iota_{B,T,U,\psi}=\Ind_{U^F}^{G^F}(\psi\circ\iota_{B,T,U,\psi})=\Ind_{U^F}^{G^F}(\psi^{-1})=\Gamma_{\psi^{-1}}$.
    This proves the Lemma.
\end{proof}

\subsection{Deligne-Lusztig theory and Jordan decomposition}\label{sec: DL theory and JD}
In this section,  we briefly recall Deligne-Lusztig theory and Lusztig's Jordan decomposition of characters of finite group of Lie type.
For more details, we refer to \cite{DLtheory},\cite{Jordandecopo}, \cite{Jordandecodisconnectedcase},\cite{Book:DigneandMichel} and \cite{Book:GeckandMalle}.

Let $\fq$ be a finite field of characteristic $p$.
Let $G$ be as before a connected reductive group over $\fq$ with the corresponding Frobenius morphism $F:G\rightarrow G$.
\begin{Fact}[{\cite[Proposition 11.1.16]{Book:DigneandMichel}}]\label{remark geometric conjugacy}
Let $T$ be a maximal $F$-stable torus of $G$ and let $\theta$ be a linear character of $T^F$.
The $G^F$-conjugacy classes of pairs $(T,\theta)$ are in one-to-one correspondence
with the ${G^{*}}^{F^{*}}$-conjugacy classes of pairs $(T^*, s)$ where $s$ is a semisimple
element of $G^{*^{F^*}}$ and $T^*$ is an $F^*$-stable maximal torus of $G^*$ containing $s$.
\end{Fact}
Let $(T,\theta)$ be a pair as above.
For every such pair, Deligne–Lusztig  constructed a virtual representation of $G^F$ using $\ell$-adic cohomology with compact support ($\ell\neq p$).
This virtual representation is denoted by $R_{T}^{G}(\theta)$.
If $(T,\theta)$ and $(T^{*}, s)$ corresponds to each other as in the Fact \ref{remark geometric conjugacy} then sometimes we will write $R_{T^*}^{G}(s)$ for $R_{T}^{G}(\theta)$.
Let $L$ be a $F$-stable Levi subgroup of $G$ and $\chi$ be an  representation of $L^F$. Lusztig extended the notion of Deligne-Lusztig induction and constructed a virtual representation $R_L^G(\chi)$ of $G^F$. The functor $R_L^G$ is called the Lusztig induction functor.
If $L$ is an $F$-stable Levi subgroup of $F$-stable parabolic subgroup $P$ then $R_L^G(\chi)$ is simply the Harish-Chandra induction and in this case $R_L^G(\chi)$ is a representation of $G^F$.
\begin{definition}[{\cite[Def. 2.6.1]{Book:GeckandMalle}}]
Let $s\in {G^{*}}^{F^*}$ be semisimple. 
Define $\mathcal{E}(G^F,s)$ to be the set of all $\rho\in\Irr(G^F)$ such that 
$\langle R_{T^*}^G(s), \rho \rangle\neq 0$ for some $F^*$-stable maximal torus $T^{*}\subseteq G^{*}$ with $s\in T^*$.
This set is called a \textit{rational series} of characters of $G^F$, or \textit{Lusztig series} of characters.
\end{definition}
\begin{definition}[Unipotent character]
A character $\rho\in\Irr(G^F)$ is called a \textit{unipotent}
character if $\langle R_{T^*}^G(1),\rho\rangle\neq 0$ for some $F^*$-stable maximal torus $T^*\subseteq G^*$.
We will denote the set of unipotent characters of $G^F$ by $\Uch(G^F)$ and by definition we have $\E(G^F,1)=\Uch(G^F)$.
\end{definition}
\begin{Fact}[{\cite[Prop.11.1.1, Prop.11.3.2]{Book:DigneandMichel}}]
For any $\rho\in\Irr(G^F)$, there exists an $F$-stable maximal
torus $T$ and $\theta\in\Irr(T^F)$ such that $\langle\rho,R_{T}^{G}(\theta)\rangle\neq 0$.
Moreover we have,
$$\Irr(G^F)=\coprod_{(s)}\E(G^F,s),$$
where $(s)$ runs over
the semisimple conjugacy classes of $G^{*^{F^*}}$.
\end{Fact}
\begin{definition}[{\cite[Definition 7.1.5]{Book:DigneandMichel}}]
Let $F$ be a Frobenius morphism on the connected reductive group $G$.
Define $\epsilon_G$ to be the sign $(-1)^{F\mathrm{-rank}(G)}$, where $F\mathrm{-rank}(G)$ denotes the $F$-split rank of $G$.
\end{definition}
In \cite{Jordandecopo}, Lusztig introduced the notion of Jordan decomposition of characters for a connected reductive group with connected centre.
This gives us a parameterization of $\mathcal{E}(G^F,s)$ by the unipotent characters of $C_{G^*}(s)^{F^*}$.
\begin{theorem}[Lusztig]\label{thm:Lusztig}
Assume that $Z(G)$ is connected.
Let $s\in G^{*^{F^*}}$ be a semisimple element, and let $H\coloneqq C_{G^*}(s)$.  There is a bijection,
$$\E(G^F,s)\rightarrow \Uch(H^{F^*}),\ \ \ \ \rho\mapsto u_{\rho},$$
such that for any $F^*$ stable maximal torus  $T^*\subseteq H$, we have,
$$\langle R_{T^*}^{G}(s),\rho\rangle=\epsilon_G\epsilon_H\langle R_{T^*}^{H}(1_{T^*}),u_\rho\rangle.$$
\end{theorem}
A collection of bijections $\{\E(G^F,s)\rightarrow \Uch(H^{F^*})\}_{(s)}$, satisfying Theorem \ref{thm:Lusztig}  is called a  \textit{Jordan decomposition of characters}.
Jordan decomposition of characters is not unique in general. It leaves some ambiguity.
This ambiguity was resolved by Digne and Michel who proved the following uniqueness result: 
\begin{theorem}[{\cite[Theorem 7.1]{lusztigsparametrazitationdigne}}]\label{Jordandecompositionconnectedcentre}
There exists a unique collection of bijections,
$$J_{s}^{G}:\mathcal{E}(G^F,s)\rightarrow \Uch(C_{G^*}(s)^{F^*}),$$
where $G$ runs over connected reductive groups with connected centre and Frobenius
map $F$, and $s\in G^{*^{F^*}}$ is semisimple, satisfying the following, where we write $H \coloneqq C_{G^*}(s)$:
\begin{enumerate}
    \item  For any $F^*$-stable maximal torus $T^*\leq H$,
    $$\langle R_{T^*}^{G}(s),\rho\rangle=\epsilon_G\epsilon_H\langle R_{T^*}^{H}(1_{T^*}),J_{s}^{G}(\rho)\rangle,\ \ \ \ \ \ \mbox{for all}\ \rho\in\mathcal{E}(G^F,s).$$
    \item If $s=1$ and $\rho\in\Uch(G^F)$ is unipotent then,
    \begin{enumerate}
        \item[a)] the Frobenius eigenvalues $\omega_{\rho}$ and $\omega_{J_{1}^{G}(\rho)}$ are equal, and 
        \item[b)] if $\rho$ lies in the principal series then $\rho$ and $J_{1}^{G}(\rho)$ correspond to the same character of the Iwahori-Hecke algebra.
    \end{enumerate}
    \item If $z\in Z(G^{*^{F^*}})$ then $J_{sz}^{G}(\rho\otimes \hat{z})=J_{s}^{G}(\rho)$ for $\rho\in\mathcal{E}(G^F, s)$, where $\hat{z}$ is the linear character of $G^F$ corresponding to $z$.
    \item For any $F^*$-stable Levi subgroup $L^*$ of $G^*$ such that $H\leq L^*$, with dual $L\leq G$, the following diagram commutes:
  \[ \begin{tikzcd}
\mathcal{E}(G^F,s)\arrow{r}{J_{s}^G}
& \Uch(H^{F^*})  \\%
\mathcal{E}(L^F,s) \arrow{u}{R_{L}^{G}}\arrow{r}{J_{s}^L}
& \Uch(H^{F^*})\arrow{u}{\Id}.
\end{tikzcd}
\]
\item If $G$ is of type $E_8$ and $H$ is of type $E_7\times A_1$ (resp. $E_6\times A_2$) and $L\leq G$ is a Levi subgroup of type $E_7$ (resp. $E_6$) with dual $L^*\leq H$ then the following diagram
commutes:
\[ \begin{tikzcd}
\mathbb{Z}\mathcal{E}(G^F,s)\arrow{r}{J_{s}^G}
&\mathbb{Z} \Uch(H^{F^*})  \\
\mathbb{Z}\mathcal{E}(L^F,s)_{c} \arrow{u}{R_{L}^{G}}\arrow{r}{J_{s}^L}
& \mathbb{Z}\Uch(L^{*^{F^*}})_{c}\arrow{u}{R_{L^*}^{H}},
\end{tikzcd}
\]
where the index $c$ denotes the subspace spanned by the cuspidal part of the
corresponding Lusztig series.
\item For any $F$-stable central torus $T_{1}\leq Z(G)$ with corresponding natural epimorphism $\pi_1:G\rightarrow G_{1}\coloneqq G/T_{1}$ and for $s_1\in G^{*^{F^*}}$
 with $s =\pi_1^*(s_1)$ the following
diagram commutes:
 \[ \begin{tikzcd}
\mathcal{E}(G^F,s)\arrow{r}{J_{s}^G}
& \Uch(H^{F^*}) \arrow{d}{} \\%
\mathcal{E}(G_{1}^F,s_1) \arrow{u}{}\arrow{r}{J_{s_{1}}^{G_1}}
& \Uch(H_{1}^{F^*}),
\end{tikzcd}
\]
with $H_1=C_{G_{1}^{*}}(s_{1})$ and where the vertical maps are just the inflation map along $G^F\rightarrow G_{1}^F$ and the restriction along the embedding $H_{1}^{F^*}\rightarrow H^{F^*}$ respectively.
\item If $G$ is a direct product $\prod_{i}G_{i}$ of F-stable subgroups $G_i$ then $J_{\prod s_{i}}^{G}=\prod J_{s_{i}}^{G_{i}}$.
\end{enumerate}

\end{theorem}

Later in \cite{Jordandecodisconnectedcase}, Lusztig extended this  Jordan decomposition of characters to the general case (\textit{i.e.},  allowing $Z(G)$ to be disconnected).
The main ingredient of his proof is  reduction to the case where $Z(G)$ is connected using a regular embedding.
We will now state some standard facts.
For more details we refer to \cite{Jordandecodisconnectedcase}.

\begin{Fact}
Let $\pi:G\rightarrow G_{\ad}$ be the adjoint quotient of $G$.
We have a natural isomorphism $G_{\ad}^F/\pi(G^F)\cong (Z(G)/Z(G)_{\circ})_{F}$
(the subscript $F$ denotes $F$-coinvariants, i.e. largest quotient on which $F$ acts trivially).
The group $G_{\ad}^F$ acts naturally on $G^F$ by automorphisms:
$$g:g_1\mapsto \dot{g}g_{1}\dot{g}^{-1},\ \ \mbox{where}\ \pi(\dot{g})=g.$$ 
This induces a natural action of $G_{\ad}^F$  on $\Irr(G^F)$. 
Note that under this action,  $\pi(G^F)$ acts trivially $\Irr(G^F)$.
Therefore, we have an action of $G_{\ad}^F/G^F$ on $\Irr(G^F)$ (recall our notation \ref{sec:notation}).
This action can be extended by linearity to virtual representations.
One can easily show that this extended action stabilizes each $R_{T^*}^G(s)$.
Thus, we have an action of $G_{\ad}^F/G^F$ $\left(\cong (Z(G)/Z(G)_{\circ})_F\right)$ on $\E(G^F,s)$.
\end{Fact}
\begin{definition}[regular embedding]
Let $G$ be a connected reductive group over $\fq$.
Then a morphism $i:G\rightarrow G'$ is called a regular embedding if
$G'$ is a connected reductive group over $\fq$ with connected centre, $i$ is an isomorphism of $G$ with a closed subgroup of $G'$ and $i(G)$, $G'$ have the same derived subgroup.

\end{definition}

Let $i:G\rightarrow G'$ be a regular embedding.
By \cite[Theorem 1.7.12]{Book:GeckandMalle}, it corresponds to a surjective homomorphism $i^*:G'^*\rightarrow G^*$ (over $\fq$).
Let us denote the kernel of $i^*$ by $K$, which is a central torus in $G'^*$.
We have a natural isomorphism: $$K^{F^*}\xrightarrow{\cong}\Hom(G'^F/G^F,\Qlcl^{\times}).$$
This induces an action of $K^{F^*}$ on $\Irr(G'^F)$ by tensor product.
Under this action, $k\in K^{F^*}$ maps $\E(G'^F,s')$ to $\E(G'^F,ks')$.
Define $K_{s'}^F$ to be the  set of all $k\in K$ which map $\E(G'^F,s')$ into itself or, equivalently,
$$K_{s'}^{F^*}=\{k\in K^{F^*}: ks'\ \mbox{is conjugate to}\ s'\ \mbox{under}\ G'^{*^{F^*}} \}.$$
\begin{Fact}\label{rmk action Ks'}
Let $s'\in G'^{*^{F^*}}$ such that $s=i^*(s')\in G^{*^{F^*}}$.
Denote the centralizer of $s$ in $G^*$ by $H$ as before, \textit{i.e.}, $H\coloneqq C_{G^*}(s)$.
Then there is a natural isomorphism $H^{F^*}/H^{F^*}_{\circ}\cong K^{F^*}_{s'}$ given by the correspondence $x\in H^{F^*} \mapsto s'^{-1}\dot{x}s'\dot{x}^{-1}\in K^{F^*},$ where $\dot{x}\in G'^*$ such that $i^*(\dot{x})=x$.
Using this isomorphism, we have an action of  $H^{F^*}/H^{F^*}_{\circ}$ on $\E(G'^{*^{F^*}}, s')$.
\end{Fact}

\begin{Fact}[{\label{rmk; action HF}\cite[Prop.2.3.15]{Book:GeckandMalle}}]
Let $i:G\rightarrow G'$ be a regular embedding  as above and let $s'\in G'^{*^{F^*}}$ such that $s=i^*(s')\in G^{*^{F^*}}$.
Let $H'$ be the centralizer of $s'$ in $G'^*$.
Then $i^*$ defines a surjective homomorphism from $H'$ onto $H_{\circ}$ with kernel $K$.
 Hence, we have a canonical bijection 
$$\Uch(H_{\circ}^{F^*})\rightarrow\Uch(H'^{F^*}),\ \ \ \rho'\mapsto\rho'\circ i^*|_{H'^{F^*}}.$$
Using this, the action (by conjugation) of $H^{F^*}/H^{F^*}_{\circ}$ on $\Uch(H_{\circ}^{F^*})$
becomes an action of $H^{F^*}/H^{F^*}_{\circ}$ on $\Uch(H_{\circ}'^{F^*})$.
\end{Fact}
\begin{lemma}[{\cite[Prop.8.1]{Jordandecodisconnectedcase}}]\label{rmk: compatible action}
Let $s'\in G'^{*^{F^*}}$ be semisimple and let $i^*(s')=s\in G^{*^{F^*}}$. 
Let $H'=C_{G'^*}(s')$ and $H=C_{G^*}(s)$.
Then the unique bijection,
 $$J_{s'}^{G'}:\E(G'^F,s')\rightarrow\Uch(H'^{F^*}),$$
is compatible with the action of $H^{F^*}/H_{\circ}^{F^*}$.
\end{lemma}
\begin{proof}
This follows from property $3$ of $J_s'$ in the Theorem \ref{Jordandecompositionconnectedcentre} and by the definition of action on both sides.
\end{proof}
\begin{Fact}\label{rmk action G'}
Let $i:G\rightarrow G'$ be a regular embedding. 
Then there is a natural surjective homomorphism ${G'^{F}}/G^F\rightarrow G_{\ad}^F/G^F$.
Thus, $G'^{F}/G^F$ acts on $\E(G^F,s)$ via $G_{\ad}^F/G^F$.
\end{Fact}

\begin{lemma}[{\cite[\textsection 10]{Jordandecodisconnectedcase}}]
Let $i:G\rightarrow G'$ be a regular embedding.
Then for any $\rho'\in \Irr(G'^F)$, the restriction $\rho'|_{G^F}$ is multiplicity free.
\end{lemma}
Next, we  state a result from representation theory of finite groups:
\begin{Fact}[{\cite[\textsection 9]{Jordandecodisconnectedcase}}]\label{rmk:action multi}
Let $N$ be a normal subgroup of a finite group $H$ such that $H/N$ is abelian.
Then there is a natural action of the abelian group $H/N$ on $\Irr(N)$ and 
there is a natural action of the abelian group $\widehat{H/N}$ on $\Irr(H)$, where $\widehat{H/N}$ is the Pontryagin dual of $H/N$.

Assume that any $\rho'\in\Irr(H)$ restricts to a multiplicity free representation of $N$. 
Then there is a unique bijection:
$$\Irr(N)/\sim{H/N}\longleftrightarrow \Irr(H)/\sim{\widehat{H/N}},$$
with the following properties. 
Let $\mathcal{O}$ be an $H/N$-orbit on $\Irr(N)$ and let $\mathcal{O}'$ be the
corresponding $\widehat{H/N}$-orbit on $\Irr(H)$.
If $\rho_{\circ}'\in \mathcal{O}'$ and $\rho_{\circ}\in\mathcal{O}$
then 
$$\rho_{\circ}'|_{N}=\sum_{\rho\in \mathcal{O}}\rho\ \ \mbox{and}\ \ \Ind_{N}^{H}(\rho_{\circ})=\sum_{\rho'\in\mathcal{O}'}\rho'.$$
Moreover, the stabilizer of $\rho_{\circ}$ in $H/N$ and the
stabilizer of $\rho_{\circ}'$ in $\widehat{H/N}$ are orthogonal to each other under the natural duality $H/N \times \widehat{H/N}\rightarrow\Qlcl^{\times}$.
\end{Fact}
With these notations, we will state the result of Lusztig.

\begin{theorem}[{\cite[Proposition 5.1 ]{Jordandecodisconnectedcase}}]\label{Thm: JD disconnected}
There exists a surjective map,
$$J_s=J_s^{(G,G')}:\E(G^F,s)\rightarrow \Uch(H_{\circ}^{F^*})/\sim{H^{F^*}/H_{\circ}^{F^*}},$$ with the following properties:
\begin{enumerate}
    \item The fibres of $J_s$ are precisely the orbits of the action of $G_{\ad}^F/G^F$ on $\E(G^F,s)$.
\item If $\mathcal{O}$ is an $H^{F^*}/H_{\circ}^{F^*}$-orbit on $\Uch(H_{\circ}^{F^*})$  and $\Gamma\leq H^{F^*}/H_{\circ}^{F^*}$ is the stabilizer of an element in $\mathcal{O}$, then the fibre $J_s^{-1}(\mathcal{O})$ has
precisely $|\Gamma|$ elements. 
\item If $\rho\in J_s^{-1}(\mathcal{O})$ and $T^*$ is an $F^*$-stable maximal torus of $G^*$ containing $s$, then
$$\langle R_{T^*}^{G}(s),\rho\rangle_{G}=\epsilon_G\epsilon_{H_{\circ}}\sum_{\rho'\in \mathcal{O}}\langle R_{T^*}^{H}(1_{T^*}),\rho'\rangle_{H_{\circ}}.$$
\end{enumerate}
\end{theorem}

The map $J_s$ is defined as follows:
\begin{center}
 \begin{tikzcd}
 \E(G^F,s)\arrow[d, "\mbox{Fact}\ \ref{rmk action G'}"]\\
 \E(G^F,s)/\sim\ {G'}^F/G^F\arrow[d, "\mbox{Fact}\ \ref{rmk:action multi}"]\\
 \bigcup_{k\in K^{F^*}}\E(G'^F, s'k)/\sim K^{F^*}\cong\widehat{({G'}^F/G^F)}\arrow[d, "\mbox{Fact}\ \ref{rmk action Ks'}"]\\
 \E(G'^F, s')/\sim K^{F^*}_{s'}=H^{F^*}/H_{\circ}^{F^*}\arrow[d, "\mathrm{Lemma}\ \ref{rmk: compatible action}"]\\
 \Uch(H'^{F^*})/\sim (H^{F^*}/H_{\circ}^{F^*})\arrow[d, "\mbox{Fact}\ \ref{rmk; action HF}"]\\
 \Uch(H_{\circ}^{F^*})/ \sim (H^{F^*}/H_{\circ}^{F^*})
\end{tikzcd}
\end{center} 
\begin{Fact}[{\cite[Theorem 12.4.12]{Book:DigneandMichel}}]\label{factgelfandchar}
For a given algebraic Whittaker datum $(B,\psi)$ and any semi-simple conjugacy class $(s)$
of $G^{*^{F^*}}$,  there is exactly one $\psi$-generic representation $\gamma_{s,\psi}$ in $\E(G^F,s)$ and we have,
$$\Gamma_{\psi}=\sum_{(s)}\gamma_{s,\psi},$$
where $(s)$ varies over the set of semisimple classes of $G^{*^{F^*}}$.
\end{Fact}

\subsection{Parameterization of unipotent representation }
In this section we will recall the Lusztig's parameterization of  unipotent representations
in terms of the set $\mathfrak{X}(W,F)$.
By standard reduction argument one can reduce the classification of unipotent representations to the case where $G$ is a simple of adjoint type. 
For more details about this reduction argument we refer to \cite[Remark 4.2.1]{Book:GeckandMalle}.

Before proceeding further, let us recall the definition of \textit{Frobenius eigenvalue} associated to a unipotent representation.
For more details we refer to \cite[\textsection 4.2.21]{Book:GeckandMalle}.
Let $\delta$ be a smallest integer such that $F^{\delta}$ acts trivially on $W$. 
For $w\in W$, $X_w$ denotes the Deligne-Lusztig variety associated to $w$.
\begin{Fact}[{\cite[Corollary 3.9]{Lusztig77}},{\cite[\textsection 4.2.21]{Book:GeckandMalle}}] \label{thm:eigenvalue}
Let $G$ be a connected reductive group defined over $\fq$.
Assume that $G$ is simple of adjoint type.
Let $\rho$ be an  unipotent representation of $G^F$.
Then it is clear that there exist $w\in W$, $i\geq 0$ and $\mu\in\Qlcl^{\times}$ such that  $\rho$ is isomorphic to the $G^F$-submodule of $H^i_c(X_w,\Ql)_{\mu}$ (the generalised $\mu$-eigenspace of $F^{\delta}$ on $H^i_c(X_w,\Ql)\otimes\Qlcl$).
Suppose that $\rho$ is also isomorphic to the $G^F$-submodule of $H^i_c(X_{w'},\Ql)_{\mu'}$ for some $\mu'$.
Then $\mu'\cdot\mu^{-1}$ is an integral power of $q^{\delta}$.
Thus, $\mu$ is uniquely determined by $\rho$ upto a factor which is an integral power of  $q^{\delta}$.   

Furthermore, there is a well-defined root of unity $\omega_{\rho}$
such that $$\mu=\omega_{\rho}q^{m\delta/2}$$ 
where $m\in\mathbb{Z}$. 
(Here, we assume that a square root $q^{1/2}\in\Qlcl$ has been fixed.) 
\end{Fact}
\begin{definition}[Frobenius eigenvalue]
Let $\rho$ be a unipotent representation of $G^F$. Let $\omega_{\rho}$ be the root of unity associated to $\rho$ as in Fact \ref{thm:eigenvalue}.
We call $\omega_{\rho}$ the \textit{Frobenius eigenvalue} of $\rho$.
\end{definition}

Now we briefly recall the description of a set $\mathfrak{X}(W,F)$.
This description is in two steps.
In the first step, we describe a set
$\mathfrak{X}^{\circ}(W,F)\subset\mathbb{C}^{\times}\times\mathbb{Z}$, which is only dependent on the pair $(W,F)$.
If $W$ is trivial, then $\mathfrak{X}^{\circ}(W,F)=\{(1,1)\}$.
This set $\mathfrak{X}^{\circ}(W,F)$ parameterizes the cuspidal unipotent representations of $G^F$.

\begin{Fact}[\cite{Jordandecopo},{\cite[Theorem 4.5]{Geck17}}]\label{para:cuspidal}
Let $G$ be a connected reductive group defined over $\fq$.
Assume that $G$ is simple of adjoint type.
Then there exists a unique bijection $ \Uch(G^F)_c\rightarrow\mathfrak{X}^{\circ}(W,F)$  with the following property.
If $\rho\in\Uch(G^F)_c$ corresponds to $x=(\omega,m)\in\mathfrak{X}^{\circ}(W,F)$, then $\omega=\omega_{\rho}$ and $m=n_{\rho}$, where $\omega_\rho$ is the Frobenius eigenvalue of $\rho$ and $n_{\rho}$ is the smallest  positive integer such that $n_{\rho}\mathbb{D}_{\rho}\in \mathbb{Z}[\mathbf{q}]$.
Here $\mathbb{D}_{\rho}$ denotes the \textit{degree polynomial} (see \cite[Definition 2.3.25]{Book:GeckandMalle}) in the indeterminate $\mathbf{q}$. 
\end{Fact}

Let $S\subset W$ be the set of simple reflections determined $B_0$.
Let $\mathcal{P}(S)$  be the set of all subsets of $S$.
Then $F$ acts on $\mathcal{P}(S)$ and we denote by $\mathcal{P}(S)^F$
the $F$-stable subsets.
For $J\in\mathcal{P}(S)^F$, define $W_{J}=\{w\in W^F:wJw^{-1}=J\}$.
We now define the set $\mathfrak{X}(W,F)$ to be, 
$$\mathfrak{X}(W,F)\coloneqq\{(J,\phi,x):J\in\mathcal{P}(S)^F,\phi\in\Irr(W_J), x\in\mathfrak{X}^{\circ}(W_J,F) \}.$$

Let $J\in\mathcal{P}(S)^F$ and let $P$ be the $F$-stable parabolic subgroup determined by $J$ with Levi factor $L$.
Let $\tau$ be a cuspidal unipotent representation of $L^F$.
Define $W_{L,\tau}\coloneqq N_{W_{G}}(W_L)^F/W_L^F$. Then we have $W_J\cong W_{L,\tau}$.
We will denote the Harish-Chandra series associated to a cuspidal pair $(L,\tau)$ by $\Irr(G^F,(L^F,\tau))$.
By \cite[Theorem 3.2.5]{Book:GeckandMalle}, there is  a natural bijection $\Irr(G^F,(L^F,\tau))\xrightarrow{\sim}\Irr(W_{L,\tau})$), $\rho\mapsto\tau_{\rho}$.
\begin{Fact}
Let $G$ be connected reductive group defined over $\fq$.
Assume that $G$ is simple of adjoint type.
Then there exists a unique bijection $\Uch(G^F)\rightarrow\mathfrak{X}(W,F)$ with the following property.
If $\rho\in\Uch(G^F)$ corresponds to $(J,\phi,x)\in\mathfrak{X}(W,F)$ then $\phi=\tau_{\rho}$, where $\tau\in \Uch(L_J^F)_c$ corresponds to $x$ under the bijection in Fact \ref{para:cuspidal}.
\end{Fact}

\section{Jordan decomposition and Automorphism}
Let $G$ be a connected reductive group defined over a finite field $\fq$ and let $F:G\rightarrow G$ be the corresponding Frobenius morphism.
Let $G^*$ be the dual reductive group defined over $\fq$ equipped with Frobenius morphism $F^*$ determined by the choice of a pinning $\p$ of $G$.  Then $(G,F)$ and $(G^*, F^*)$ are in duality.
Let $\sigma:G\rightarrow G$ be an $F$ automorphism of $G$ that respects $\p$ and let $\sigma^*:G^*\rightarrow G^*$ be the corresponding dual $F^*$-automorphism.
Then $\sigma$ (resp. $\sigma^*$) induces an automorphism of the finite group $G^F$ (resp. $G^{*^{F^*}}$).
This induces an action of $\sigma$ on $\Irr(G^F)$ as $\rho\mapsto\rho\circ \sigma^{-1}$.
Then, $$R_{T}^{G}(\theta)\circ \sigma^{-1}=R_{\sigma(T)}^{G}(\theta\circ \sigma^{-1}).$$
Let $\rho\in \E(G^F,s)$. For any pair $(T,\theta)$ which is in duality with $(T^*,s)$ (recall Fact \ref{remark geometric conjugacy}), we have,  
$$0\neq \langle \rho, R_{T}^{G}(\theta)\rangle =\langle  \rho\circ \sigma^{-1}, R_{T}^{G}(\theta)\circ \sigma^{-1}\rangle=\langle  \rho\circ \sigma^{-1}, R_{\sigma(T)}^{G}(\theta\circ \sigma^{-1})\rangle.$$
Note that $\sigma$ maps the $G^F$-conjugacy class of a pair $(T,\theta)$ to the $G^F$-conjugacy class of a pair $(\sigma(T),\theta\circ \sigma^{-1})$.
Let $(T^*,s)$ be a pair corresponding to the pair $(T,\theta)$ as in Fact \ref{remark geometric conjugacy}.
By \cite[Prop. 7.2]{actionofauto}, $(\sigma^{*^{-1}}(T^*), \sigma^{*^{-1}}(s))$ corresponds to the pair $(\sigma(T),\theta\circ \sigma^{-1})$ and
$\sigma$ induces a bijection,
$$f_{\sigma}:\mathcal{E}(G^F,s)\rightarrow \mathcal{E}(G^F,\sigma^{*^{-1}}(s))\ \ \ \rho\mapsto \rho\circ \sigma^{-1}.$$
Similarly, $\sigma^{*^{-1}}$ induces a bijection,
$$f_{\sigma^{*}}^{-1}: \Uch (C_{G^*}(\sigma^{*^{-1}}(s)))\rightarrow\Uch(C_{G^*}(s));\ \ \ u_{\rho}\mapsto u_{\rho}\circ \sigma^{*^{-1}}.$$

The following Theorem proves the finite group counterpart of the first half of the  \cite[Conjecture 1]{dualityinvolution}. 
\begin{theorem}\label{thm:generic}
Let $G$ be a connected reductive reductive group defined over $\fq$.
Let $\gamma_{s,\psi}$ be a  $\psi$-generic representation of $G^F$.
Then, 
$$\gamma_{s,\psi}\circ\iota_{G,\p}=\gamma_{s,\psi}^{\vee},$$
where $\p$ is the pinning corresponding to the algebraic Whittaker datum $\psi$.
\end{theorem}
\begin{proof}
Note that if $\gamma_{s,\psi}\in \E(G^F,s)$, then $\gamma_{s,\psi}\circ\iota_{G,\p},\ \gamma_{s,\psi}^{\vee}\in\E(G^F,s^{-1})$.
By Lemma \ref{lemma generic}, $\gamma_{s,\psi}\circ\iota_{G,\p}$ is $\psi^{-1}$-generic and $\gamma_{s,\psi}^{\vee}$ is also $\psi^{-1}$-generic.
Hence by Fact \ref{factgelfandchar},  we get, $$\gamma_{s,\psi}\circ\iota_{G,\p}=\gamma_{s^{-1},\psi^{-1}}=\gamma_{s,\psi}^{\vee}.$$
This proves the result. 
\end{proof}
\begin{corollary}\label{cor:generic}
    Let $G$ be a connected reductive group defined over $\fq$.
    Let $\rho_{s,\psi}$ be a semisimple character of $G^F$.
    Then $$\rho_{s,\psi}\circ\iota_{G,\p}=\rho_{s^{-1},\psi^{-1}}=\rho_{s,\psi}^{\vee},$$ where $\p$ is the pinning corresponding to the algebraic Whittaker datum $\psi$. 
\end{corollary}
\begin{proof}
    The result follows from Theorem \ref{thm:generic} and properties of the operator $D_G$.
\end{proof}

Let  $H\coloneqq C_{G^*}(s)$, the centralizer of $s$ in $G^*$.
Let us define a new bijection,
$$\mathbb{J}_{s}^{G}:\E(G^F, s)\rightarrow\Uch(H^{F^*})\ \  \mbox{as}\ \  \J_{s}^{G}(\rho)=f_{\sigma^{*}}^{-1}\circ J_{\sigma^{*^{-1}}(s)}^G\circ f_{\sigma}(\rho) =J_{\sigma^{*^{-1}}(s)}^G(\rho\circ \sigma^{-1})\circ \sigma^{*^{-1}}.$$
We have the following result:
\begin{theorem}\label{thm:connectedcentre}
Let $G$ be a connected reductive group over $\fq$ such that $Z(G)$ is connected.
Then, $$\J_{s}^{G}=J_{s}^{G},$$
where $J_s^G$ is the unique Jordan decomposition as in the  Theorem
\ref{Jordandecompositionconnectedcentre}.
\end{theorem}
\begin{proof}
We will show that the new bijection $\J_s^{G}$ satisfies the properties $(1)-(7)$   in the Theorem \ref{Jordandecompositionconnectedcentre}.
\begin{enumerate}
    \item Let $\rho\in\E(G^F,s)$ be any element.
    We have, for any $F^*$-stable maximal torus $T^*\leq H$, \begin{align*}
        \langle R_{T^*}^{H}(1_{T^*}),\J_{s}^G(\rho)\rangle_{H}
        &=\langle R_{T^*}^{H}(1_{T^*}),J_{\sigma^{*^{-1}}(s)}^G(\rho\circ \sigma^{-1})\circ \sigma^{*^{-1}}\rangle_{H}\\
         &=\langle R_{T^*}^{H}(1_{T^*})\circ \sigma^{*} ,J_{\sigma^{*^{-1}}(s)}(\rho\circ \sigma^{-1})\rangle_{C_{G^*}(\sigma^{*^{-1}}(s))}\\
         &=\langle R_{\sigma^{*^{-1}}(T^{*})}^{C_{G^*}(\sigma^{*^{-1}}(s))}(1_{\sigma^{*^{-1}}(T^{*})}),J_{\sigma^{*^{-1}}(s)}(\rho\circ \sigma^{*^{-1}})\rangle_{C_{G^*}(\sigma^{*^{-1}}(s))}\\
         &=\epsilon_{G}\epsilon_{C_{G^*}(\sigma^{*^{-1}}(s))}\langle R_{\sigma^{*^{-1}}(T^{*})}^{G}(\sigma^{*^{-1}}(s)),\rho\circ \sigma^{-1}\rangle_{G}\\
         &=\epsilon_{G}\epsilon_{H}\langle R_{T^{*}}^{G}(s),\rho)\rangle_{G}.
         \end{align*}
         Thus, we get $\langle R_{T^{*}}^{G}(s),\rho)\rangle_{G}=\epsilon_{G}\epsilon_{H}\langle R_{T^*}^{H}(1_{T^*}),\J_{s}^G(\rho)\rangle_{H}$ for any $\rho\in\E(G^F,s)$. 
         \item As $\sigma$ and $\sigma^*$ commute with the Frobenius morphism, so
          the Frobenius eigenvalues corresponding to $\rho$ (resp. $J_1(\rho)$) and $\rho\circ \sigma^{-1}$ (resp. $J_{1}(\rho)\circ \sigma^{*^{-1}}$) are equal. This follows since  $\sigma$ (resp. $\sigma^*$) induces an isomorphism between the cohomology groups associated with Deligne-Lusztig varieties which commutes with the action of Frobenius on cohomology.
         By Theorem \ref{Jordandecompositionconnectedcentre}, the Frobenius eigenvalues corresponding to $\rho\circ \sigma^{-1}$ and $J_{1}(\rho\circ \sigma^{-1})$
         are equal. Hence the Frobenius eigenvalues corresponding to $\rho$ and $\J_1(\rho)=J_{1}(\rho\circ \sigma^{-1})\circ \sigma^{*^{-1}}$ are equal.
         This proves the first part of $(2)$.
         
         For the second part of $(2)$, we assume $\rho$ is in the principal series, meaning that it is a constituent of $\Ind_{B^F}^{G^F}(1)$, and so $J_1(\rho)$ is also in the principal series, and they correspond to the same character $\chi$ of the Hecke algebra $H(G^F, B^F)$ (which is identified with  $H(G^{*^{F^*}}, B^{*^{F^*}}) $ via the natural isomorphism $\delta:W_G(T_0)\rightarrow W_{G^*}(T_0^*)$ of Weyl groups).
         If $\rho$ (resp. $J_1(\rho)$) is in the principal series, then so is $\rho\circ \sigma^{-1}$ (resp. $J_1(\rho)\circ \sigma^{*^{-1}}$).
         By Theorem \ref{Jordandecompositionconnectedcentre}, $\rho\circ \sigma^{-1}$ and $J_1(\rho\circ \sigma^{-1})$ corresponds to same character of the Iwahori-Hecke algebra.
         Note that $\sigma$ induces an isomorphism between the Hecke algebra $H(G^F, B^F)$ and $H(G^F, \sigma(B)^F)$. 
         Under this isomorphism, the character corresponding to $\rho$ goes to the character corresponding to $\rho\circ \sigma^{-1}$.
         Similarly, $\sigma^*$ induces an isomorphism between the Hecke algebra $H(G^{*^{F^*}}, B^{*^{F^*}}) $ and $H(G^{*^{F^*}}, \sigma^*(B^{*^{F^*}}))$. 
         Under this isomorphism the character corresponding $J_1(\rho)$ goes to the character corresponding to $J_1(\rho)\circ \sigma^{*^{-1}}$.
         Thus $\rho$ and $\J_1(\rho)=J_1(\rho\circ \sigma^{-1})\circ \sigma^{*^{-1}}$ corresponds to the same character of the Iwahori-Hecke algebra $H(G^F, B^F)$ (which is identified with  $H(G^{*^{F^*}}, B^{*^{F^*}}) $ through the natural isomorphism $\delta:W_G(T_0)\rightarrow W_{G^*}(T_0^*)$ of Weyl groups).
         
\item This part  follows from the fact that, if $\hat{z}$ is a character corresponding to $z\in Z(G^{{*}^{F^*}})$, then $\hat{z}\circ \sigma^{-1}$ is a character corresponding to $\sigma^{*^{-1}}(z)$.
\item Let  $L^*$ be an $F$-stable Levi subgroup of $G^*$ such that $H\leq L^*$, with dual $L\leq G$.
Then by Theorem \ref{Jordandecompositionconnectedcentre} we have $J_{s}^L=J_{s}^{G}\circ R_{L}^{G}$.
Note that $\sigma$ maps $\E(L^F,s)$ to $\E(\sigma(L)^F,\sigma^{*^{-1}}(s))$ and also for $\rho\in\E(L^F,s)$, we have $R_{\sigma(L)}^{G}(\rho\circ\sigma^{-1})=R_{L}^{G}(\rho)\circ\sigma^{-1}$.
We have, \begin{align*}
    \J_{s}^G\circ R_{L}^{G}
    &=f_{\sigma^{*}}^{-1}\circ J_{\sigma^{*^{-1}}(s)}^L\circ f_{\sigma}\circ R_{L}^{G}\\
    &=f_{\sigma^{*}}^{-1}\circ J_{\sigma^{*^{-1}}(s)}^G\circ R_{\sigma(L)}^{G}\circ f_{\sigma}\\
    &=f_{\sigma^{*}}^{-1}\circ J_{\sigma^{*^{-1}}(s)}^{\sigma^*(L)}\circ f_{\sigma}\\
    &= \J_{s}^{L}.
\end{align*}
 This implies that the following diagram commutes:
  \[ \begin{tikzcd}
\mathcal{E}(G^F,s)\arrow{r}{\J_{s}^G}
& \Uch(H^{F^{*}})  \\%
\mathcal{E}(L^F,s) \arrow{u}{R_{L}^{G}}\arrow{r}{\J_{s}^L}
& \Uch(H^{F^{*}})\arrow{u}{\Id}.
\end{tikzcd}
\]
\item Suppose $G$ is of type $E_8$ and $H$ is of type $E_7 A_1$ (resp. $E_6 A_2$) and $L\leq G$ is a Levi subgroup of type $E_7$ (resp. $E_6$) with dual $L^*\leq H$.
Then by Theorem \ref{Jordandecompositionconnectedcentre},
$$J_{s}^{G}\circ R_{L}^{G}(\rho)=R_{L^*}^{H}\circ J_{s}^L(\rho),\ \ \ \mbox{for any}\ \rho\in \mathbb{Z}\mathcal{E}(L^F,s)_{c}$$
where the index $c$ denotes the subspace spanned by the cuspidal part of the corresponding Lusztig series.
We have,
\begin{align*}
    \J_{s}^G\circ R_{L}^{G}(\rho)
    &=f_{\sigma^{*}}^{-1}\circ J_{\sigma^{*^{-1}}(s)}^G\circ f_{\sigma}\circ R_{L}^{G}(\rho)\\
    &=f_{\sigma^{*}}^{-1}\circ J_{\sigma^{*^{-1}}(s)}^G\circ R_{\sigma(L)}^{G}\left(f_{\sigma}(\rho)\right)\\
   &=f_{\sigma^{*}}^{-1}\circ R_{\sigma^*(L^*)}^{H}\circ J_{\sigma^{*^{-1}}(s)}^{\sigma(L)} \left(f_{\sigma}(\rho)\right)\\
   &= R_{L^*}^{H}\circ f_{\sigma^{*}}^{-1}\circ J_{\sigma^{*^{-1}}(s)}^{\sigma(L)}\circ f_{\sigma}(\rho)\\
   &=R_{L^*}^{H}\circ \J_{s}^L (\rho).
\end{align*}
This implies the following diagram commutes:
\[ \begin{tikzcd}
\mathbb{Z}\mathcal{E}(G^F,s)\arrow{r}{\J_{s}^G}
&\mathbb{Z} \Uch(H^{F^*})  \\%
\mathbb{Z}\mathcal{E}(L^F,s)_{c} \arrow{u}{R_{L}^{G}}\arrow{r}{\J_{s}^L}
& \mathbb{Z}\Uch(L^{*^{F^*}})_{c}\arrow{u}{R_{L^*}^{H}}.
\end{tikzcd}
\]
\item Let $T\leq Z(G)$ be any $F$-stable central torus and let $\pi_1:G\rightarrow G_{1}\coloneqq G/T$ be the  natural epimorphism.
Let $T_1=\sigma(T)$ be the image of $T$ under $\sigma$, so $T_1\leq Z(G)$ is also an $F$-stable central torus in $G$.
Let $\pi_2:G\rightarrow G_{2}\coloneqq G/T_{1}$ be the  corresponding natural epimorphism.
Then $\sigma^{-1}$ induces a natural morphism $f:\Irr(G_{1}^F)\rightarrow\Irr(G_{2}^F)$.
For $s_i\in G_i^*$ 
 with $s =\pi_1^*(s_1)$ and $\sigma^{*^{-1}}(s)=\pi_{2}^{*}(s_2)$, by Theorem \ref{Jordandecompositionconnectedcentre} the following
diagrams commute:
\begin{center}
     \begin{tikzcd}
\mathcal{E}(G^F,s)\arrow{r}{J_{s}^G}
& \Uch(H^{F^*}) \arrow{d}{\psi_1} \\%
\mathcal{E}(G_{1}^F,s_1) \arrow{u}{\phi_1}\arrow{r}{J_{s_{1}}^{G_1}}
& \Uch(H_{G_1}^{F^*})
\end{tikzcd}\ \ \ \
 \begin{tikzcd}
\mathcal{E}(G^F,\sigma^{*^{-1}}(s))\arrow{r}{J_{\sigma^{*^{-1}}(s)}^G}
& \Uch(H^{F^*}) \arrow{d}{\psi_2} \\%
\mathcal{E}(G_{2}^F,s_2) \arrow{u}{\phi_2}\arrow{r}{J_{s_{2}}^{G_2}}
& \Uch(H_{G_2}^{F^*}),
\end{tikzcd}
 \end{center}
 
where $H_{G_i}=C_{G_{i}^{*}}(s_{i})$ for $i=1,2$ and where the vertical maps are just the inflation map along $G^F\rightarrow G_{1}^F$ (resp.$G^F\rightarrow G_{2}^F$) and the restriction along the embedding $H_{G_1}^{F^*}\rightarrow H^{F^*}$ (resp. $H_{G_2}^{F^*}\rightarrow H^{F^*}$) respectively.
There is a natural bijection $f^*:\Uch(H_{G_1}^{F^*})\rightarrow\Uch(H_{G_2}^{F^*})$, induced by $\sigma^*$.
Then it follows from the definition of $f_{\sigma}$ and $f$ that,
$$\phi_2\circ f=f_{\sigma}\circ \phi_1\ \ \ \mbox{and}\ \ \ {f_{\sigma^*}}^{-1}\circ\psi_2=\psi_1\circ f_{\sigma^*}^{-1} .$$
We have,
\begin{align*}
    \psi_1\circ\J_s^G\circ \phi_1
    &= \psi_1\circ f_{\sigma^{*}}^{-1}\circ J_{\sigma^{*^{-1}}(s)}^L\circ f_{w}\circ\phi_1\\
     &= f_{\sigma^*}^{-1}\circ\psi_2\circ J_{\sigma^{*^{-1}}(s)}^G\circ\phi_2\circ f\\
     &= f_{\sigma^*}^{-1}\circ J_{s_2}^{G_2}\circ f\\
     &= \J_{s_1}^{G_1}.
\end{align*}
This implies the following diagram commutes:
\[
 \begin{tikzcd}
\mathcal{E}(G^F,s)\arrow{r}{\J_{s}^G}
& \Uch(H^{F^*}) \arrow{d}{\psi_1} \\%
\mathcal{E}(G_{1}^F,s_1) \arrow{u}{\phi_1}\arrow{r}{\J_{s_{1}}^{G_1}}
& \Uch(H_{G_1}^{F^*}).
\end{tikzcd}
\]
\item For the final property $(7)$,
suppose $s=\prod_is_i$.
Consider \begin{align*}
    \J_s^G &= f_{\sigma^{*}}^{-1}\circ J_{\sigma^{*^{-1}}(s)}^G\circ f_{\sigma}\\
    &=f_{\sigma^*}^{-1}\circ \prod_i J_{\sigma^{*^{-1}}(s_i)}^{\sigma(G_i)}\circ f_{\sigma}\\
     &=\prod_i f_{\sigma^{*}}^{-1}\circ J_{\sigma^{*^{-1}}(s_i)}^{\sigma(G_i)}\circ f_{\sigma}\\
     &=\prod_i \J_{s_i}^{G_i}.
\end{align*}
Hence $(7)$ follows.
\end{enumerate}
Thus, the new bijection $\J_{s}^G$ also satisfies the properties satisfied by $J_s$.
By Theorem \ref{Jordandecompositionconnectedcentre}, we get $\J_s^G=J_s^G$.
This proves the lemma.
\end{proof}

\begin{lemma}\label{lem:unipiota}
Let $G$ be a connected reductive group over $\fq$ and let $\rho\in\Uch({G^F})$. Then $\rho\circ\iota_{G,\p}$ is independent of $\p$.
\end{lemma}

\begin{proof}
This follows from the fact that $\Uch(G^F)$ and $\Uch({G_{\ad}}^F)$ are in natural bijection and $\iota_{{G_{\ad},\p}}$ is independent of $\p$.
\end{proof}

\begin{theorem}\label{thm: comm of involution }
Let $G$ be a connected reductive group over $\fq$ with connected centre.
Fix a pinning $\p=(G,B,T,\{X_{\alpha}\})$ of $G$ and let $\iota_{G,\p}$ be the corresponding duality involution (which does not  depend on the choice of pinning).
Then for any $\rho\in\E(G^F,s)$,
$$J_{\iota_{G^*,\p^*}(s)^{}}(\rho\circ \iota_{G,\p}^{})=J_{s}(\rho)\circ \iota_{G^*,\p^*}.$$
\end{theorem}
\begin{proof}
The Theorem follows from Theorem \ref{thm:connectedcentre} and 
Fact \ref{rmk:dual of duality}.
\end{proof}
\begin{theorem}\label{Prop commofdual}
    Let $G$ be a connected reductive group over $\fq$ with connected centre.
    Let $s\in G^{*^{F^*}}$ be a semisimple element. 
    Then for any $\rho\in \E(G^F,s)$, 
    $$J_{s^{-1}}(\rho^{\vee})=J_s(\rho)^{\vee},$$
    where $J_s$ is the unique Jordan decomposition as in Theorem \ref{Jordandecompositionconnectedcentre}.
\end{theorem}
\begin{proof}
Define a new bijection, $$\J_{s}^G:\E(G^F,s)\rightarrow\Uch(H^{F^*}) \ \ \mbox{as}\ \ \J_{s}^{G}(\rho)=J_{s^{-1}}^G(\rho^{\vee})^{\vee}.$$
Then by arguments analogous to those in Theorem \ref{thm:connectedcentre}, one can easily show that the new bijection $\J_{s}^G$ satisfies the  properties $(1)-(7)$ of Theorem \ref{Jordandecompositionconnectedcentre}.
The result follows.
\end{proof}


\subsection{General case}
Let $G$ be a connected reductive group over $\fq$ and let $Z(G)$ denote the centre of $G$.
Let $i:G\hookrightarrow G'$ be a regular embedding constructed as follows.
Let $S\subset G$ be an $F$-stable torus such that $Z(G)\subseteq S$.
Let $G'$ be the quotient of $G\times S$ by the closed normal subgroup
$\{(z,z^{-1}):z\in Z(G)\}$.
Let $S'$ be the image of $\{1\}\times S\subseteq G\times S$ in $G'$.
Then the map $i:G\rightarrow G'$ induced by $G\rightarrow G\times S$, $g\mapsto (g,1)$ is a regular embedding and $S'$ is the centre of $G'$.
Let $T\subset G$ be an $F$-stable maximal torus, then $T'=T\cdot S'$ is an $F$-stable maximal torus of $G'$.
Also if $B$ is an $F$-stable Borel subgroup of $G$ containing $T$ then $B'=B\cdot S'$ is an $F$-stable Borel subgroup of $G'$ containing $T'$.
Thus $(G,B_0,T_0)$ determines a triple $(G',B'_0,T'_0)$.
Moreover, if $\p=(G,B,T,\{X_{\alpha}\})$ is a pinning of $G$ then $\p'=(G',B',T',\{X_{\alpha}\})$ is a pinning of $G'$.
From the construction, we can observe that the following diagram commutes:
\begin{center}
\begin{tikzcd}
G\arrow[r,"c_{G,\p}"]\arrow[hookrightarrow]{d}[swap]{i}
& G\arrow[hookrightarrow,"i"]{d}\\
G'\arrow[r,"c_{G',\p'}"]
& G'.
\end{tikzcd}
\end{center}
Let $\{J_s^G=J_s^{(G,G')}\}$ be the collection of surjections as in the Theorem \ref{Thm: JD disconnected}.

\begin{theorem}\label{Thm invocommdisconn}
Let $G$ be a connected reductive group over $\fq$.
Assume that $G$ satisfies the condition in Prop. \ref{prop pinnindependence}.
Let us fix a pinning $\p=(G,B,T,\{X_{\alpha}\})$ of $G$ and let $\iota_{G,\p}$ be the corresponding duality involution (which does not depend on the choice of our pinning).
Then for any $\rho\in\E(G^F,s)$,
$$J^{G}_{s^{-1}}(\rho\circ \iota_{G,\p})=J^{G}_{s}(\rho)\circ \iota_{G^*,\p^*}.$$
\end{theorem}
\begin{proof}
Since the duality involution takes a semisimple element to a conjugate of its inverse,
 it maps $\E(G,s)$ to $\E(G,s^{-1}).$
Note that the duality involution is independent of the pinning and $G_{\ad}^F$ acts transitively on the set of pinnings.
The action of $G'^F$ on $\Irr(G^F)$ is via $G_{\ad}^F$, hence it induces a map, 
$$\E(G^F,s)/\sim ({G'}^F/G^F)\rightarrow\E(G^F,s)/\sim ({G'}^F/G^F)\ \ \mbox{given by}\ \ [\rho]\mapsto[\rho\circ\iota_{G,\p}].$$
Let $i:G\rightarrow G'$ be a regular embedding and $i^*:G'^*\rightarrow G^*$ be the corresponding surjective morphism with kernel $K$.
Since $i^*\circ \iota_{G'^*,\p'*}=\iota_{G^*,\p^*}\circ i^*$, hence $\iota_{G'^*,\p'^*}$ maps $K$ to $K$.
Let $s'\in G'^*$ such that $i^*(s')=s$.
Then $\iota_{G',\p'}$ maps  $\E(G'^F, s'k)/\sim K^{F^*}$ to $\E(G'^F, (s'k)^{-1})/\sim K^{F^*}_{s'}$ as $[\rho]\mapsto[\rho\circ\iota_{G',\p'}]$.
Also, $\iota_{G',\p'}$ maps  $\E(G'^F, s')/\sim K^{F^*}$ to $\E(G'^F, (s')^{-1})/\sim K^{F^*}_{s'}$ as $[\rho]\mapsto[\rho\circ\iota_{G',\p'}]$.

Similarly,  $\iota_{G'^*,\p'^*}$ maps   $\Uch(H^{F^*}_{G'})/\sim (H^{F^*}/(H)_{\circ}^{F^*})$ to $\Uch(H^{F^*}_{G'})/\sim (H^{F^*}/(H)_{\circ}^{F^*})$. (By abuse of notations we denote the restriction of $\iota_{G'^*,\p'^*}$ to $H^{F^*}_{G'}$ by $\iota_{G'^*,\p'^*}$ ).
To prove the Theorem it is enough to show that each square in the following diagram commutes:
\[ \begin{tikzcd}
\E(G^F,s)\arrow[d, "a"]\arrow[r, "\iota_{G,\p}"]
&\E(G^F,s^{-1})\arrow[d, "a"] \\
\E(G^F,s)/\sim {G'}^F/G^F\arrow[d, "b"]\arrow[r, "\iota_{G,\p}"]
& \E(G^F,s^{-1})/\sim {G'}^F/G^F\arrow[d, "b"]\\
\bigcup_{k\in K^{F^*}}\E(G'^F, s'k)/\sim K^{F^*}\cong\widehat{({G'}^F/G^F)}\arrow[d, "c"]\arrow[r, "\iota_{G',\p'}"]
&\bigcup_{k\in K^{F^*}}\E(G'^F, (s'k)^{-1})/\sim K^{F^*}\cong\widehat{({G'}^F/G^F)}\arrow[d, "c"]\\
\E(G'^F, s')/\sim K^{F^*}_{s'}=H^{F^*}/(H)_{\circ}^{F^*}\arrow[d, "d"]\arrow[r, "\iota_{G',\p'}"]
&\E(G'^F, (s')^{-1})/\sim K^{F^*}_{s'}=H^{F^*}/(H)_{\circ}^{F^*}\arrow[d, "d"]\\
\Uch(H^{F^*}_{G'})/\sim H^{F^*}/(H)_{\circ}^{F^*}\arrow[d, "e"]\arrow[r, "\iota_{G'^*,\p'^*}"]
&\Uch(H^{F^*}_{G'})/\sim H^{F^*}/(H)_{\circ}^{F^*}\arrow[d, "e"]\\
\Uch((H)_{\circ}^{F^*})/\sim H^{F^*}/(H)_{\circ}^{F^*}\arrow[r, "\iota_{G^*,\p^*}"]
&\Uch(H)_{\circ}^{F^*})/\sim H'^{F^*}/(H')_{\circ}^{F^*}.
\end{tikzcd}
\]

The first and last squares commute by definition.
Note that $\Ind_{G^F}^{G'^F}(\rho\circ\iota_{G,\p})=\Ind_{G^F}^{G'^F}(\tau)\circ\iota_{G'.\p'}$.
Then by Fact \ref{rmk:action multi}, the second square commutes.
Note that the induced action of $\iota_{G',\p'}$ on $\widehat{{G'^F}/G^F}$ under the isomorphism transfers to the action of $\iota_{G'^*,\p'^*}$ on $K^{F^*}$.
Then commutation of the third squares follows from the definition of maps.
The vertical map $d$ is the bijection $J_{s'}^{G'}$ in the case of Jordan decomposition of characters for connected reductive groups with connected center (cf. Lemma \ref{rmk: compatible action}).
By Theorem \ref{thm:connectedcentre}, $J_{s'}^{G'}$ commute with the action of $\iota_{G',\p'}$ and $\iota_{G'^*,\p'^*}$ respectively.
Hence the square involving $d$  commutes.

This implies that the surjective map in the Theorem \ref{Thm: JD disconnected} commutes with the corresponding action of $\iota_{G',\p'}$ and $\iota_{G'^*,\p'^*}$.
This proves the Theorem.
\end{proof}

\begin{theorem}\label{thm:commofdualdiscon}
    Let $G$ be a connected reductive group over $\fq$.
    Let $s\in G^{*^{F^*}}$ be a semisimple element. 
    Then for any $\rho\in \E(G^F,s)$, 
    $$J_{s^{-1}}(\rho^{\vee})=J_s(\rho)^{\vee},$$
    where $J_s$ is the Jordan decomposition as in Theorem \ref{Thm: JD disconnected}.
\end{theorem}
\begin{proof}
    This follows from Theorem \ref{Prop commofdual} and by an argument similar to the proof of Theorem \ref{Thm invocommdisconn}.
\end{proof}

\section{Main Results}
We continue with the notations of the previous section. 
\begin{lemma}\label{lem:iota}
    Let $G$ be a connected reductive group defined over $\fq$ with connected centre.
    Let $\iota_{G}$ be the duality involution. Then
    for any $F$-stable torus $T$, the pair $(\iota_{G}(T),\theta\circ\iota_{G})$ is $G^F$-conjugate to $(T,\theta^{-1})$
\end{lemma}
\begin{proof}
In the view of the Fact \ref{remark geometric conjugacy} and  \cite[Prop. 7.2]{actionofauto}, it is enough to show that, for any $F$-stable torus $T$ and a semisimple element $t\in T^F$,  the pair $(\iota_{G}(T),\iota_{G}(t))$ is $G^{F}$-conjugate to $(T,t^{-1})$.
Fix a pinning $\p=(G,B_0,T_0,{X_{\alpha}})$. 
Let $w\in W$ such that $T$ is of type $w$. Fix a representative $\dot{w}$ of $w$ in $N_G(T_0)$ and a representative $w_G$ of the longest Weyl group element in $N_G(T_0)^F$ as in \cite[\S 2.2]{Changyang}. It then follows from \cite[Theorem 3.4]{Changyang} that $\iota_{G}(\dot{w})=w_G\dot{w}w_G^{-1}$.
By Lang–Steinberg theorem, there exists a $g\in G$ such that $\dot{w}=g^{-1}F(g)$ and 
$(T,t)=\prescript{g}{}{}(T_0,t_0)=(gT_0g^{-1},gt_0g^{-1})$ for some $t_0\in T_0$. 
Then, we have,
\begin{align*}
(\iota_{G}(T), \iota_{G}(t))
& =(\iota_{G}(gT_0g^{-1}),\iota_{G}(gt_0g^{-1}))\\
& =(\iota_{G}(g)T_0\iota_{G}(g)^{-1},\iota_{G}(g)w_Gt_0^{-1}w_G^{-1}\iota_{G}(g)^{-1})\\
& =\prescript{g_0}{}{}(T,t^{-1})
\end{align*}
where $g_0=\iota_{G}(g)w_Gg^{-1}$.
We know that, $\iota_{G}(g^{-1}F(g))=w_Gg^{-1}F(g)w_G^{-1}$.
This imples that,
 $F(\iota_{G}(g)w_Gg^{-1})=\iota_{G}(g)w_Gg^{-1}$.
 Hence,  $g_0=\iota_{G}(g)w_Gg^{-1}\in G^F$. 
This proves the Lemma.
\end{proof}

\begin{theorem}
Let $G$ be a connected reductive group over $\fq$ with connected centre.
Fix a pinning $\p=(G,B,T,\{X_{\alpha}\})$ of $G$ and let $\iota_{G,\p}$ be the corresponding duality involution (which does not depend on the choice of pinning).
Then for any $\rho\in \E(G^F,s)$,
$$\rho\circ\iota_{G,\p}=\rho^{\vee}\ \ \mbox{if and only if}\ \ J_s(\rho)\circ\iota_{G^*,\p^*}=J_{s}(\rho)^{\vee}.$$
\end{theorem}
\begin{proof}
The proof follows from Theorem \ref{thm: comm of involution } and Theorem \ref{Prop commofdual}.
\end{proof}
 \begin{lemma}\label{lemma:stabilizer}
 Let $G$ be a connected reductive group over $\fq$ and let $F:G\rightarrow G$ be a corresponding Frobenius morphism.
 For a cuspidal irreducible representation $\rho$ of $G^F$ there exists a semisimple character $\rho_{s}$ in the rational Lusztig series of $\rho$ having the
same stabilizer as $\rho$ in the group $G_{\ad}^F/G^F$.
 \end{lemma}
 \begin{proof}
 \begin{enumerate}
     \item Assume that $G$ is simple and simply connected.
     Then the lemma follows from \cite[Theorem 1]{mallecuspodalcharacter}.
     \item Next, assume that $G$ is a semisimple simply connected type.
     Then $G\cong\prod_i G_{i}$, where $G_i$ is a simple simply connected type and $G_{\ad}=\prod_i(G_{i})_{\ad}$.
     Suppose, $\rho$ is an irreducible cuspidal representation of $G$ then $\rho=\otimes_i\rho_i$, where $\rho_i$ is an irreducible cuspidal representation of $G_i$.
     Then the Lemma follows  by applying $(a)$ to each $G_i$.
     \item Assume that $G$ is a connected reductive group over $\fq$ such that $G_{\der}$ is semisimple of simply connected type.
     Let $\rho\in \Irr(G^F)$ be a cuspidal irreducible representation of $G^F$.
     Then $\rho|_{G_{\der}^F}=\sum_{x\in G^F/G_{\der}^F}\prescript{x}{}{}\chi=\sum_{x\in G^F/G_{\der}^F}\chi$, since for a coset of $G_{\der}^F$ in $G^F$ we can choose a coset representative in $Z(G)$.
     Note that $G_{\ad}=(G_{\der})_{\ad}$ and there is a natural surjective morphism from $(G_{\der})_{\ad}^F/G_{\der}^F$ to $G_{\ad}^F/G^F$ whose kernel is $G^F/G_{\der}^F$.
     Hence, we have an action $G_{\ad}^F/G^F$ on $\Irr(G^F_{\der})$ via $(G_{\der})_{\ad}^F/G_{\der}^F$.
     Then, $\Stab_{G_{\ad}^F/G^F}(\chi)=\Stab_{(G_{\der})_{\ad}^F/G_{\der}^F}(\chi)$.
     By part $(b)$, there exists a semisimple character $\chi_s$ in the same rational Lusztig  series of $\chi$ such that $\Stab_{(G_{\der})_{\ad}^F/G_{\der}^F}(\chi)=\Stab_{(G_{\der})_{\ad}^F/G_{\der}^F}(\chi_s)$.
     Let $\rho_s$ be a semisimple character in the rational Lusztig series of $\rho$ such that $\rho_s|_{G_{\der}^F}=\sum_{x\in G^F/G_{\der}^F}\chi_{s}$.
     
     Let $x\in \Stab_{G_{\ad}^F/G^F}(\rho)$ be any element. 
     Then $x$-stabilizes the rational Lusztig series of $\rho$ and $x\in\Stab_{G_{\ad}^F/G^F}(\chi)$.
     Therefore, $x$ stabilises the semisimple character $\chi_s$ (since $\Stab_{G_{\ad}^F/G^F}(\chi)=\Stab_{G_{\ad}^F/G^F}(\chi_s)$).
     If $x$ does not fix $\rho_s$, then $\prescript{x}{}{}\rho_s$ and $\rho_s$  differ by a character of $G^F/G_{\der}^F$.
     This contradicts the fact that $x$-stabilizes the rational Lusztig series of $\rho$.
     Hence, $x\in\Stab_{G_{\ad}^F/G^F}(\rho_s)$ and $\Stab_{G_{\ad}^F/G^F}(\rho)\subseteq\Stab_{G_{\ad}^F/G^F}(\rho_s)$.
     Similarly, we can show that $\Stab_{G_{\ad}^F/G^F}(\rho_s)\subseteq\Stab_{G_{\ad}^F/G^F}(\rho)$.
     This implies \newline
     $\Stab_{G_{\ad}^F/G^F}(\rho_s)=\Stab_{G_{\ad}^F/G^F}(\rho)$.
      \item By \cite[1.7.13]{Book:GeckandMalle}, there exists a connected reductive group $G'$ over $\fq$ and a surjective morphism $f:G'\rightarrow G$  defined over $\fq$ such that $\ker(f)$ central torus, $G_{\der}$ is semisimple simply connected type and $f$ induces a surjective homomorphism $f:G'^F\rightarrow G^F$. 
     Hence, we can identify the irreducible representations of $G^F$ with the irreducible representations of $G'^F$ which are trivial on $\ker(f)^F$.
     Also, there is a induced surjective homomorphism from $G'^F_{\ad}/G'^F\rightarrow G_{\ad}^F/G^F$ whose kernel contained in centre.
     This induces an action of $G_{\ad}^F/G^F$ on $\Irr(G'^F)$ via $G'^F_{\ad}/G'^F$.
     Also, if $\rho'_s$ is a semisimple character of $G'^F$, then there exists a semisimple character $\rho_s$ of $G^F$ such that $\rho'_s=\rho_s\circ f$.
     Let $\rho$ be an irreducible cuspidal representation of $G^F$ and let $\rho'=\rho\circ f$. 
     Then $\Stab_{G_{\ad}^F/G^F}(\rho)=\Stab_{G_{\ad}^F/G^F}
     (\rho')=\Stab_{G'^F_{\ad}/G'^F}(\rho')$.
     The result then follows from part $(c)$.
     \end{enumerate}
     This proves the Lemma.
 \end{proof}
 
 \begin{lemma}\label{lemma:addualiota}
Let $G$ be a connected reductive group over $\fq$ with corresponding Frobenius $F$ and let $s$ be a semisimple element of $G^{*^{F^*}}$.
Assume that $s=1$ or $2\mathrm{H}^1(F,Z(G))=0$.
Fix a pinning $\p=(G,B,T,\{X_{\alpha}\})$ of $G$ and let $\iota_{G,\p}$ denote the corresponding duality involution (which does not depend on the pinning chosen if $s\neq1$).
Let $\rho\in\E(G^F,s)$ be any element such that $\rho\circ\iota_{G,\p}=\rho^{\vee}\circ\ad(g)$ for some $g\in G_{\ad}^F$.
Then $\rho\circ\iota_{G,\p}=\rho^{\vee}$.
\end{lemma}
\begin{proof}

We will first prove the result for irreducible cuspidal representations. 
Let $\rho\in\E(G^F,s)$ be a cuspidal representation of $G^F$. Then by Lemma \ref{lemma:stabilizer}, there exists a semisimple character $\rho_{s,\psi}\in\E(G^F,s)$ such that 
 $$\Stab_{G_{\ad}^F/G^F}(\rho)=\Stab_{G_{\ad}^F/G^F}(\rho_{s,\psi}).$$

For a semisimple element $s\in {G^*}^{F^*}$, define $A(s)\coloneqq C_{G^*}(s)/C_{G^*}(s)_{\circ}$ and $\widehat{A(s)^{F^*}}$ to be the Pontryagin dual of $A(s)^{F^*}$.
Recall there is a  surjective homomorphism  $\zeta:G_{\ad}^F/G^F\rightarrow \widehat{A(s)^{F^*}}$ such that $\ker(\zeta)=\Stab_{G_{\ad}^F/G^F}(\rho_{s,\psi})$ (see \cite[Prop 3.12]{Digne92}).
This implies that $\mathcal{O}_{\rho_{s,\psi}}$ is a torsor over $\widehat{A(s)^{F^*}}$.  
Fix $\Gad$-equivariant bijections $\gamma_{1,\psi}:\AS\rightarrow\mathcal{O}_{\rho_{s,\psi}}$
and $\gamma_{2,\psi}:\AS\rightarrow\mathcal{O}_{\rho_{s,\psi}}$ determined
by the maps $1\mapsto\rho_{s,\psi}$ and $1\mapsto\rho_{s,\psi}^{\vee}\circ\iota$
respectively. We have $\Gad$-equivariant bijections $\eta_{1}:\mathcal{O}_{\rho}\rightarrow\mathcal{O}_{\rho_{s,\psi}}$
and $\eta_{2}:\mathcal{O}_{\rho}\rightarrow\mathcal{O}_{\rho_{s,\psi}}$
determined by $\rho\mapsto\rho_{s,\psi}$ and $\rho^{\vee}\circ\iota\mapsto\rho_{s,\psi}^{\vee}\circ\iota_{G,\p}$
respectively. These bijections induce $\Gad$-equivariant bijections
$\gamma_{1}:\AS\rightarrow\mathcal{O}_{\rho}$ and $\gamma_{2}:\AS\rightarrow\mathcal{O}_{\rho}$
determined by $1\mapsto\rho$ and $1\mapsto\rho^{\vee}\circ\iota_{G,\p}$
respectively. 

Now consider the bijection $D:\tau\in\mathcal{E}(G,s)\mapsto\tau^{\vee}\circ\iota_{G,\p}\in\mathcal{E}(G,s)$.
This bijection is $\Gad$-equivariant and consequently induces $\Gad$-equivariant
bijections $D_{\rho_{,\psi}}:\mathcal{O}_{\rho,\psi}\rightarrow\mathcal{O}_{\rho,\psi}$
and $D_{\rho}:\mathcal{O}_{\rho}\rightarrow\mathcal{O}_{\rho}$. Let
$\bar{D}_{\rho,\psi}:\AS\rightarrow\AS$ and $\bar{D}_{\rho}:\AS\rightarrow\AS$
be the bijections induced by $D_{\rho,\psi}$ and $D_{\rho}$ respectively.
The former are then necessarily $\Gad$-equivariant. 

We claim that $\bar{D}_{\rho,\psi}=\bar{D}_{\rho}=\mathrm{id}$ .
From the commutative diagram:
\[
\xymatrix{\mathcal{O}_{\rho}\ar@{->}[rrr]^{D_{\rho}}\ar@{->}[dd]^{\eta_{1}} & \ar@{->}[r] & \ar@{-}[r] & \mathcal{O}_{\rho}\ar@{->}[dd]^{\eta_{2}}\\
 & \AS\ar@{->}[lu]^{\gamma_{1}}\ar@{->}[r]_{\bar{D}_{\rho,\psi}}^{\bar{D}_{\rho}}\ar@{->}[ld]^{\gamma_{1,\psi}} & \AS\ar@{->}[ru]^{\gamma_{2}}\ar@{->}[rd]^{\gamma_{2,\psi}}\\
\mathcal{O}_{\rho_{s,\psi}}\ar@{->}[rrr]^{D_{\rho,\psi}} &  &  & \mathcal{O}_{\rho_{s,\psi}}
},
\]
we get that $\bar{D}_{\rho,\psi}=\bar{D}_{\rho}$ . From Corollary \ref{cor:generic}, $D_{\rho,\psi}=\mathrm{id}$. We thus get that $\bar{D}_{\rho,\psi}=\bar{D}_{\rho}=\mathrm{id}$.  If  $\rho\circ\iota_{G,\p}=\rho^{\vee}\circ\ad(g)$ for some $g\in G_{\ad}^F$, then $\bar{D}_{\rho}:1\mapsto \zeta(g^{-1})$. But since $\bar{D}_{\rho}=\mathrm{id}$, $g^{-1}\in \mathrm{ker}(\zeta)=\Stab_{G_{\ad}^F/G^F}(\rho)$. Thus   $D_{\rho}=\mathrm{id}$. This completes the proof in the case when $\rho$  is cuspidal. 

Now let $\rho\in \E(G^F,s)$ to be an arbitrary element.
Then there exists a cuspidal pair $(L,\tau)$, $\tau\in\E(L^F,s)$ such that $\rho$ belongs to $\Irr(G^F,(L^F,\tau))$, the \textit{Harish-Chandra series} associated to $(L,\tau)$, (refer to \cite[Defintion 3.1.14]{Book:GeckandMalle} for the definition of Harish-Chandra serires). Since $H^1(F,Z(G))$ is an elementary abelian $2$-group by hypothesis, so is $H^1(F,Z(L))$ by \cite[Lemma 12.3.5]{Book:DigneandMichel}.

By the first part, we have $\tau\circ\iota_{L,\p}=\tau^{\vee}$.
Also, $R_{L}^{G}(\tau\circ\iota_{L,\p})=R_{L}^{G}(\tau)\circ\iota_{G,\p}=R_{L}^{G}(\tau^{\vee})$.
Hence, $\rho\circ\iota_{G,\p}$ and $\rho^{\vee}$  belong to the same Harish-Chandra series.
By hypothesis, $\rho^{\vee}$ and $\rho\circ\iota_{G}$ are $G_{\ad}^F/G^F$-conjugate.
By \cite[Corollary 3.1.17]{Book:GeckandMalle}, for them to lie in the same Harish-Chandra series, they must necessarily be equal.
This implies that, $$\rho\circ\iota_{G,\p}=\rho^{\vee}.$$
This proves the lemma.

\end{proof}


\begin{theorem}\label{thm:dualiffeigenvaluepm}
Let $G$ be a connected reductive group defined over $\fq$ with corresponding Frobenius morphism $F$.
Fix a pinning $\p=(G,B,T,\{X_{\alpha}\})$ of $G$ and let $\iota_{G,\p}$ denote the corresponding duality involution.
Let $\rho\in\Uch(G^F)$ be a unipotent representation. Then
$\rho\circ\iota_{G,\p}=\rho^{\vee}$ if and only if the Frobenius eigenvalues corresponding to $\rho$ are in $\{\pm1\}$.
\end{theorem}
\begin{proof}
By standard reduction arguments we can assume that $G$ is simple of adjoint type (see \cite[1.18]{Lusztig76} or \cite[\textsc{4.2}]{Geck17}).
Let $W$ denote the absolute Weyl group of $G$. 
Recall that there is a natural parameterization of irreducible unipotent representations in terms of the set $\mathfrak{X}(W,F)$ (see  \cite{Lusztig2015}, \cite[Corollary 4.8]{Geck17} ).
In this parameterization,  a cuspidal representation $\rho$ associates to the  pair $(\omega_\rho,n_\rho)$, where $\omega_\rho$ is the Frobenius eigenvalue of $\rho$ and $n_{\rho}$ is the smallest  positive integer such that $n_{\rho}\mathbb{D}_{\rho}\in \mathbb{Z}[\mathbf{q}]$.
Here $\mathbb{D}_{\rho}$ denotes the degree polynomial in the indeterminate $\mathbf{q}$. 
Recall that the Frobenius eigenvalues of cuspidal representations are $n^{\mbox{th}}$-roots of unity with $n=2,3,4 \mbox{~or~} 5$ (see \cite[\textsection 3.1]{Lusztig2015}, \cite[Theorem 4.5]{Geck17}).
Since the degree polynomials of $\rho,\ \rho\circ\iota_{G,\p}$ and $\rho^{\vee}$ are the same, therefore $n_{\rho\circ\iota_{G,\p}}=n_\rho=n_{{\rho}^{\vee}}$.
Also, we have $\omega_{\rho}=\omega_{\rho\circ\iota_{G,\p}}=\omega_{\rho^{\vee}}^{-1}$.
Thus, for cuspidal representations, $\rho\circ\iota_{G,\p}=\rho^{\vee}$ if and only if $\omega_{\rho}\in\{\pm 1\}$.

For a cuspidal pair $(L,\tau)$, define $W_{L,\tau}\coloneqq N_{W_{G}}(W_L)^F/W_L^F$. This is a Coxeter group.
Since $\iota_{G,\p}$ acts on $W_G$ by a conjugation action of the longest Weyl group element $w_0$, the induced map on $W_{L,\tau}$ is conjugation by the image of $w_0$.
Hence the action of $\iota_{G.\p}$ on $\Irr(W_{L,\tau})$ is trivial.
Since $W_{L,\tau}$ is a coxeter group, all representations of it are self-dual.
Let $\rho\in \Uch(G^F)$ be any element.
Then there exists a cuspidal pair $(L,\tau)$, $\tau\in\Uch(L^F)$, such that $\rho\in\Irr(G^F,(L^F,\tau))$ (Harish-Chandra series associated with $(L,\tau)$).
Suppose, $\rho$ corresponds to a triple $(J,\phi,(\omega_{\tau},n_{\tau}))\in\mathfrak{X}(W_G,F)$, where $J$ is the set of simple roots  corresponding to $L$ and $\phi$ is an irreducible representation of $W_{L,\tau}$ corresponding to $\rho$ (By \cite[Theorem 3.2.5]{Book:GeckandMalle}, there is a natural bijection between $\Irr(G^F,(L^F,\tau))$ and $\Irr(W_{L,\tau})$).
Then $\rho^{\vee}$ corresponds to a triple $(J,\phi^{\vee},(\omega_{\tau^{\vee}},n_{\tau^{\vee}}))=(J,\phi,(\omega_{\tau}^{-1},n_{\tau}))$ and $\rho\circ\iota_{G,\p}$ corresponds to a triple $(J,\phi\circ\iota_{G,\p},(\omega_{\tau},n_{\tau}))=(J,\phi,(\omega_{\tau},n_{\tau}))$.
This implies that $\rho\circ\iota_{G,\p}=\rho^{\vee}$ if and only if $\omega_{\tau}\in\{\pm 1\}$.
By \cite[Prop. 4.2.23]{Book:GeckandMalle}, $\omega_{\rho}=\omega_{\tau}$.
This proves the theorem.
\end{proof}
\begin{theorem}\label{thm:iotadual}
Let $G$ be a connected reductive group over $\fq$ with corresponding Frobenius morphism $F$.
Assume that $2\mathrm{H}^1(F,Z(G))=0$.
Fix a pinning $\p=(G,B,T,\{X_{\alpha}\})$ of $G$ and let $\iota_{G,\p}$ be the corresponding duality involution (which does not depend on the choice of pinning).
Then for any $\rho\in \E(G^F,s)$,
$$\rho\circ\iota_{G,\p}=\rho^{\vee}\ \ \mbox{if and only if}\ \ 
u_\rho\circ\iota_{G^*,\p^*}=u_\rho^{\vee} \mbox{~ ~for any, and therefore all~} u_\rho\in J_s(\rho). $$
Here $J_s$ is the Jordan decomposition as in the Theorem \ref{Thm: JD disconnected}.
\end{theorem}
\begin{proof}
    Proof follows from Theorem \ref{Thm invocommdisconn}, Theorem \ref{thm:commofdualdiscon}, Lemma \ref{lemma:addualiota} and Theorem \ref{thm:dualiffeigenvaluepm}.
\end{proof}

\begin{theorem}\label{thm:main}
Let $G$ be a connected reductive group defined over $\fq$ with corresponding Frobenius morphism $F$.
Assume that $2\mathrm{H}^1(F,Z(G))=0$.
Fix a pinning $\p=(G,B,T,\{X_{\alpha}\})$ of $G$ and let $\iota_{G,\p}$ be the corresponding duality involution (which does not  depend on the choice of pinning).
Let $\rho\in\E(G^F,s)$ be any element then,
$\rho\circ\iota_{G,\p}=\rho^{\vee}$ if and only if the Frobenius eigenvalue corresponding  to any $u_{\rho}\in J_s(\rho)$ is in $\{\pm1\}$.
\end{theorem}
\begin{proof}
Suppose $\rho\circ\iota_{G,\p}=\rho^{\vee}$.
Then $J_{s^{-1}}(\rho\circ\iota_{G,\p})=J_{s^{-1}}(\rho^{\vee})=J_{s}(\rho)^{\vee}$.
Let $u_{\rho}\in J_s(\rho)$ be any element. Then $u_{\rho}\circ\iota_{G^*,\p^*}=u_{\rho}^{\vee}\circ\ad(h)$, for some $h\in H^{F^*}$.
This implies that the Frobenius eigenvalues corresponding to $u_{\rho}$ and $u_{\rho}^{\vee}$ are equal \emph{i.e.}, $\omega_{u_{\rho}^{\vee}}=\omega_{u_{\rho}}$.
Note that $\omega_{u_{\rho}^{\vee}}=\omega_{u_{\rho}}^{-1}$.
This implies that $\omega_{u_{\rho}}\in\{\pm1\}$.

Suppose for $u_{\rho}\in J_s(\rho)$, $\omega_{u_{\rho}}\in\{\pm1\}$.
Then by  Theorem \ref{thm:dualiffeigenvaluepm} and  Theorem \ref{thm:iotadual}, we get
$\rho\circ\iota_{G,\p}=\rho^{\vee}$.
This proves the Theorem.
\end{proof}

\begin{remark}\label{rem:prasad}
Given a connected reductive group $G$ defined over $\fq$, Lusztig (see \cite{Jordandecopo}) showed that $\Uch(G^F)$ gets partitioned into families and these families are parameterized by \emph{special unipotent} conjugacy classes. Families should be thought of as analogous of $L$-packets in the representation theory of reductive $p$-adic groups. If $\mathcal{F}$ is a family associated to a special unipotent $u$, then members of $\mathcal{F}$ are parameterized by pairs $(x,\psi)$, where $x\in A(u)$ and $\psi\in\Irr(Z_{A(u)}(x))$. Here $A(u)$ is the \emph{Lusztig quotient group}. If $A(u)$ an elementary abelian $2$-group (which is the case if the component group of $Z_G(u)$ is such), then all elements of the associated family $\mathcal{F}$ have Frobenius eigenvalues $\pm 1$. 
\end{remark}

In Theorem \ref{thm:main}, if $G$ has no factors of exceptional type, then the Frobenius eigenvalue of $u_\rho$ is necessarily $\pm 1$ (see \cite[Prop. 4.4.32]{Book:GeckandMalle}). We therefore have the following:
\begin{corollary}\label{thm:classical}
 Let $G$ be a connected reductive group defined over $\fq$ with corresponding Frobenius morphism $F$.
Assume that $2\mathrm{H}^1(F,Z(G))=0$. Assume further that $G$ has no factors of exceptional type. Then $\iota_G$ is a dualizing involution, i.e., $\rho\circ\iota_G=\rho^{\vee}$ for each irreducible character $\rho$ of $G^F$.   
\end{corollary}

\section{Twisted sign}\label{sec:sign}
Theorems \ref{thm:generic} and \ref{thm:main} allows us to associate a twisted sign in a natural way to most irreducible character of $G^F$ for most $G$. Twisted signs have been extensively studied in literature.  We will briefly recall the definition of this sign. 

Let $\rho$  be an irreducible character of $G^F$ and let $\iota_G$ be an involution such that we have an isomorphism, $$f:\rho\circ\iota_G\rightarrow \rho^{\vee}.$$
This induces another isomorphism, $$({f\circ\iota_G})^\vee : \rho\circ\iota_G\rightarrow \rho^{\vee}.$$ 
These two isomorphisms are easily seen to differ by a sign. 
\begin{definition}
Let $\epsilon_{\rho}\in\{\pm1\}$ be the sign associated to $\rho$ which is given by the equality: $$f=\epsilon_{\rho}(f\circ\iota_G)^\vee.$$
\end{definition}
\begin{remark}When $G=\GL_n$, then it follows from the main result of Gow \cite{Gow83} that $\epsilon_{\rho}=1$ for all irreducible characters $\rho$ of $G^F$. 
\end{remark}

\section*{Acknowledgements}
The authors are thankful to Dipendra Prasad for many helpful conversations.
The authors also thankful to Gunter Malle  who
carefully read an earlier draft of this paper and suggested many improvements.
\bibliographystyle{alpha}
\bibliography{references}

\begin{thebibliography}{MgVW87}

\bibitem[DL76]{DLtheory}
Pierre Deligne and George Lusztig.
\newblock Representations of reductive groups over finite fields.
\newblock {\em Annals of Mathematics}, 103(1):103--161, 1976.

\bibitem[DLM92]{Digne92}
Fran\c{c}ois Digne, G.~I. Lehrer, and Jean Michel.
\newblock The characters of the group of rational points of a reductive group
  with nonconnected centre.
\newblock {\em J. Reine Angew. Math.}, 425:155--192, 1992.

\bibitem[DM90]{lusztigsparametrazitationdigne}
Fran\c{c}ois Digne and Jean Michel.
\newblock On {L}usztig's parametrization of characters of finite groups of
  {L}ie type.
\newblock {\em Ast\'{e}risque}, (181-182):6, 113--156, 1990.

\bibitem[DM20]{Book:DigneandMichel}
Fran\c{c}ois Digne and Jean Michel.
\newblock {\em Representations of finite groups of {L}ie type}, volume~95 of
  {\em London Mathematical Society Student Texts}.
\newblock Cambridge University Press, Cambridge, 2020.

\bibitem[Gec17]{Geck17}
Meinolf Geck.
\newblock A first guide to the character theory of finite groups of lie type,
  2017.

\bibitem[GK72]{GelKazh72}
Israel Gel\'fand and David Ka\v{z}dan.
\newblock Representations of the group {${\rm GL}(n,K)$} where {$K$} is a local
  field.
\newblock {\em Funkcional. Anal. i Prilo\v{z}en.}, 6(4):73--74, 1972.

\bibitem[GM20]{Book:GeckandMalle}
Meinolf Geck and Gunter Malle.
\newblock {\em The character theory of finite groups of {L}ie type}, volume 187
  of {\em Cambridge Studies in Advanced Mathematics}.
\newblock Cambridge University Press, Cambridge, 2020.

\bibitem[Gow83]{Gow83}
Roderick Gow.
\newblock Properties of the characters of the finite general linear group
  related to the transpose-inverse involution.
\newblock {\em Proc. London Math. Soc. (3)}, 47(3):493--506, 1983.

\bibitem[Lus78]{Lusztig77}
George Lusztig.
\newblock {\em Representations of finite {C}hevalley groups.}
\newblock American Mathematical Society, Providence, R.I.,,, 1978.
\newblock Expository lectures from the CBMS Regional Conference held at
  Madison, Wis., August 8--12, 1977.

\bibitem[Lus84]{Jordandecopo}
George Lusztig.
\newblock {\em Characters of reductive groups over a finite field}, volume 107
  of {\em Annals of Mathematics Studies}.
\newblock Princeton University Press, Princeton, NJ, 1984.

\bibitem[Lus88]{Jordandecodisconnectedcase}
George Lusztig.
\newblock On the representations of reductive groups with disconnected centre.
\newblock {\em Ast\'{e}risque}, (168):10, 157--166, 1988.

\bibitem[Lus15]{Lusztig2015}
George Lusztig.
\newblock Restriction of a character sheaf to conjugacy classes.
\newblock {\em Bull. Math. Soc. Sci. Math. Roumanie (N.S.)},
  58(106)(3):297--309, 2015.

\bibitem[Lus77]{Lusztig76}
George Lusztig.
\newblock Coxeter orbits and eigenspaces of {F}robenius.
\newblock {\em Invent. Math.}, 38(2):101--159, 1976/77.

\bibitem[Mal17]{mallecuspodalcharacter}
Gunter Malle.
\newblock Cuspidal characters and automorphisms.
\newblock {\em Adv. Math.}, 320:887--903, 2017.

\bibitem[MgVW87]{MVW87}
Colette M\oe~glin, Marie-France Vign\'{e}ras, and Jean-Loup Waldspurger.
\newblock {\em Correspondances de {H}owe sur un corps {$p$}-adique}, volume
  1291 of {\em Lecture Notes in Mathematics}.
\newblock Springer-Verlag, Berlin, 1987.

\bibitem[Pra19]{dualityinvolution}
Dipendra Prasad.
\newblock Generalizing the {MVW} involution, and the contragredient.
\newblock {\em Trans. Amer. Math. Soc.}, 372(1):615--633, 2019.

\bibitem[RV18]{RoVin18}
Alan Roche and C.~Ryan Vinroot.
\newblock A factorization result for classical and similitude groups.
\newblock {\em Canad. Math. Bull.}, 61(1):174--190, 2018.

\bibitem[Tay18]{actionofauto}
Jay Taylor.
\newblock Action of automorphisms on irreducible characters of symplectic
  groups.
\newblock {\em J. Algebra}, 505:211--246, 2018.

\bibitem[Yan20]{Changyang}
Chang Yang.
\newblock Distinguished representations, {S}hintani base change and a finite
  field analogue of a conjecture of {P}rasad.
\newblock {\em Adv. Math.}, 366:107087, 27, 2020.

\end{thebibliography}
\end{document}